\documentclass[reqno,11pt]{amsart}

\usepackage{amsthm}
\usepackage{amssymb}
\usepackage{amsmath}
\usepackage{amsfonts}
\usepackage[utf8]{inputenc}
\usepackage[english]{babel}
\usepackage[T1]{fontenc}
\usepackage{lmodern,microtype}
\usepackage{mathpazo}
\usepackage{soul,cite}
\usepackage{cancel}

\usepackage{xcolor}
\usepackage{hyperref}

\hypersetup{%
    colorlinks,%
    linkcolor = {red!60!black},%
    citecolor = {green!60!black},%
    urlcolor = {blue!60!black}%
}

\usepackage{graphicx,psfrag}
\usepackage{subcaption}
\usepackage{tikz,nicefrac,pifont,fp}
\usetikzlibrary{arrows,backgrounds,patterns,matrix,shapes,fit,calc,shadows,%
plotmarks,arrows.meta,positioning,decorations.pathmorphing,%
decorations.pathreplacing,decorations.markings,fixedpointarithmetic}

\usetikzlibrary{external}
\tikzexternalenable
\makeatletter
\let\origsection=\section 
\def\section{\@ifstar{\origsection*}{\mysection}}
\def\mysection{\@startsection{section}{1}\z@{.7\linespacing\@plus\linespacing}{.5\linespacing}{\normalfont\scshape\centering\S}}
\makeatother

\usepackage{doi}

\usepackage{datetime}
\usepackage{lineno}
\newcommand*\patchAmsMathEnvironmentForLineno[1]{%
\expandafter\let\csname old#1\expandafter\endcsname\csname #1\endcsname
\expandafter\let\csname oldend#1\expandafter\endcsname\csname end#1\endcsname
\renewenvironment{#1}%
{\linenomath\csname old#1\endcsname}%
{\csname oldend#1\endcsname\endlinenomath}}%
\newcommand*\patchBothAmsMathEnvironmentsForLineno[1]{%
\patchAmsMathEnvironmentForLineno{#1}%
\patchAmsMathEnvironmentForLineno{#1*}}%
\AtBeginDocument{%
\patchBothAmsMathEnvironmentsForLineno{equation}%
\patchBothAmsMathEnvironmentsForLineno{align}%
\patchBothAmsMathEnvironmentsForLineno{flalign}%
\patchBothAmsMathEnvironmentsForLineno{alignat}%
\patchBothAmsMathEnvironmentsForLineno{gather}%
\patchBothAmsMathEnvironmentsForLineno{multline}%
}

\usepackage{geometry}
\newtheorem{theorem}[equation]{Theorem}
\newtheorem{lemma}[equation]{Lemma}
\newtheorem{proposition}[equation]{Proposition}

\newtheorem{claim}[equation]{Claim}
\newtheorem{question}[equation]{Question}
\newtheorem{problem}[equation]{Problem}

\theoremstyle{definition}

\newtheorem{remark}[equation]{Remark}

\numberwithin{equation}{section}

\newcommand{\oldqed}{}
\def\endofFact{\hfill$\Diamond$}
\newenvironment{claimproof}{
  \renewcommand{\oldqed}{\qed}
  \renewcommand{\qed}{\endofFact}
  \begin{proof}
  \leftskip15pt\relax
}{
  \end{proof}
  \renewcommand{\qed}{\oldqed}
}

\usepackage{xcolor}
\colorlet{tn/color/defi}{red!50!black}

\newcommand{\defi}[1]{%
  \emph{\color{red!60!black}#1}%
}

\usepackage{enumerate}
\usepackage{bookmark}
\usepackage{comment}
\usepackage{paralist}

\definecolor{darkgreen}{rgb}{.17,.5,.5} 

\newcommand{\qitem}[1]{\noindent\leavevmode\hangindent1.5\parindent%
  \noindent\hbox to1.5\parindent{#1\hss}\ignorespaces}

\def\phi{\varphi}

\newcommand{\chio}{\chi_{\mathrm{o}}}
\newcommand{\chipcf}{\chi_{\mathrm{pcf}}}

\newcommand{\cF}{\mathcal{F}}
\newcommand{\cG}{\mathcal{G}}
\newcommand{\cI}{\mathcal{I}}
\newcommand{\Choose}[2]{\genfrac(){0pt}{1}{#1}{#2}}

\usepackage{thm-restate}

\newcommand{\ci}{\alpha}
\newcommand{\cb}{\beta}
\newcommand{\cg}{\gamma}

\title{Boundedness for proper conflict-free and odd colorings}

\author[Jiménez, Knauer, Lintzmayer, Matamala, Peña, Quiroz, Sambinelli, Wakabayashi, Yu, Zamora]{
A. Jiménez
\and
K. Knauer
\and
C.N. Lintzmayer
\and
M. Matamala
\and
J.P. Peña
\and
D.A. Quiroz
\and
M. Sambinelli
\and
Y. Wakabayashi
\and
W. Yu
\and
J. Zamora
}

\address[A. Jiménez]{
  Instituto de Ingenier\'ia Matem\'atica and CIMFAV --
  Universidad de Valpara\'iso \&
Millennium Nucleus for Social Data Science (SODAS) -- Chile.
  Email: \href{mailto:andrea.jimenez@uv.cl}{andrea.jimenez@uv.cl}
}

\address[K. Knauer]{
  Departament de Matem\`atiques i Inform\`atica -- Universitat de Barcelona
  \& Centre de Recerca Matemàtica -- Spain \\
  LIS, Aix-Marseille Universit\'e, CNRS, and Universit\'e de Toulon -- France.
  Email: \href{mailto:kolja.knauer@ub.edu}{kolja.knauer@ub.edu}
}
\address[C. N. Lintzmayer, M. Sambinelli]{
  Centro de Matemática, Computação e Cognição -- 
  Universidade Federal do ABC -- 
  Brazil.  
  Emails: \href{mailto:carla.negri@ufabc.edu.br}{carla.negri@ufabc.edu.br},
  \href{mailto:m.sambinelli@ufabc.edu.br}{m.sambinelli@ufabc.edu.br}
}
\address[M. Matamala]{
  Departamento de Ingeniería Matemática and Center for Mathematical Modeling -- 
  Universidad de Chile -- 
  Chile.
  Email: \href{mailto:mar.mat.vas@dim.uchile.cl }{mar.mat.vas@dim.uchile.cl}
}
\address[J. P. Peña]{
  Departamento de Ingeniería Matemática -- 
  Universidad de Chile -- 
  Chile.
  Email: \href{mailto:juanpabloalcayaga.p21@gmail.com}{juanpabloalcayaga.p21@gmail.com}
}

\address[D. A. Quiroz]{
  Instituto de Ingenier\'ia Matem\'atica  and CIMFAV --
  Universidad de Valpara\'iso -- 
  Chile.
  Email:
  \href{mailto:daniel.quiroz@uv.cl}{daniel.quiroz@uv.cl}
}

\address[Y. Wakabayashi]{
  Instituto de Matemática e Estatística --
  Universidade de São Paulo --
  Brazil.
  Email: \href{mailto:yw@ime.usp.br}{yw@ime.usp.br}
}
\address[W. Yu]{
  Department of Mathematics --
  Zhejiang Normal University --
  China.
  Email: \href{mailto:wyu@irif.fr}{wyu@irif.fr}
} 
\address[J. Zamora]{
  Departamento de Matemáticas --
  Universidad Andres Bello --
  Chile.
  Email: \linebreak \href{mailto:josezamora@unab.cl}{josezamora@unab.cl}
}

\begin{document}


\begin{abstract}
    The \emph{proper conflict-free chromatic number}, $\chipcf(G)$, of a graph~$G$ is the least positive integer~$k$ 
    such that~$G$ has a proper $k$-coloring in which for each non-isolated vertex there is a color
    appearing exactly once among its neighbors.
    The \emph{proper odd chromatic number}, $\chio(G)$, of~$G$ is the least positive integer~$k$ such that~$G$ has a
    proper coloring in which for every non-isolated vertex there is a color appearing an odd number 
    of times among its neighbors.
    We clearly have $\chi(G) \le \chio(G) \le \chipcf(G)$.
    We say that a graph class~$\cG$ is \emph{$\chipcf$-bounded} (\emph{$\chio$-bounded}) if there is 
    a function~$f$ such that $\chipcf(G) \leq f(\chi(G))$ ($\chio(G) \leq f(\chi(G))$) for every
    $G \in \cG$.
    Caro, Petru\v{s}evski, and \v{S}krekovski (2023) asked for classes that are linearly
    $\chipcf$-bounded ($\chio$-bounded) and, as a starting point, they showed that
    every claw-free graph~$G$ satisfies $\chipcf(G) \le 2\Delta(G)+1$, which implies 
    $\chipcf(G) \le 4\chi(G)+1$.
    
    In this paper, we improve the bound for claw-free graphs to a nearly tight bound 
    by showing that such a graph~$G$ satisfies $\chipcf(G) \le \Delta(G)+6$, and even
    $\chipcf(G) \le \Delta(G)+4$ if it is a quasi-line graph.
    These results also give further evidence to a conjecture by Caro, Petru\v{s}evski, and \v{S}krekovski.
    Moreover, we show that convex-round graphs and permutation graphs 
    are linearly $\chipcf$-bounded.
    For these last two results, we prove a lemma that reduces the problem of deciding if 
    a hereditary class is linearly $\chipcf$-bounded to deciding if the bipartite graphs 
    in the class are $\chipcf$-bounded by an absolute constant.
    This lemma complements a theorem of Liu (2024) and motivates us to further study 
    boundedness in bipartite graphs.
    Among other results, we show that biconvex bipartite graphs are $\chipcf$-bounded, while convex bipartite graphs are not even $\chio$-bounded, 
    and we exhibit a class of bipartite circle graphs that is linearly $\chio$-bounded but not 
    $\chipcf$-bounded.
 \end{abstract}

\maketitle

\section{Introduction}

The frequency assignment problem has motivated much research on 2-distance colorings,
in which vertices get a different color if their distance is at most~2.
In such a coloring, each color appears at most once in each neighborhood.
A weakening of this, introduced for similar applications~\cite{1181994}, is that of
\defi{conflict-free coloring}, where for each non-isolated vertex there is a color
appearing exactly once among its neighbors.
Conflict-free colorings are not necessarily proper: indeed, three colors suffice to
conflict-free color any planar graph~\cite{3colors2017}, and, as proved 
by Pyber~\cite{Py1991}, four colors suffice to conflict-free color any line-graph.

Siting between the notion of proper coloring and 2-distance coloring are
\defi{proper conflict-free colorings}, introduced by Fabrici, Lu\v{z}ar, 
Rindo\v{s}ov\'a, and Sot\'ak~\cite{FaLuRiSo2023}.
The \defi{proper conflict-free chromatic number} of a graph~$G$ is the
smallest~$k$ such that the graph has a proper conflict-free  $k$-coloring;
it is denoted by~$\chipcf(G)$.
Fabrici \textit{et al.} showed that there are planar graphs that cannot be proper
conflict-free colored with~5 colors, and conjectured this result to be tight.
They also proved that $\chipcf(G) \le 8$ if~$G$ is planar, a bound that was 
improved by Caro, Petru\v{s}evski, and \v{S}krekovski~\cite{CaPeSk2023} 
if~$G$ has high girth or small maximum average degree (see also~\cite{ChChKwPa2023,anderson2024forbflex}). 

A variant of proper conflict-free coloring is that of \defi{proper odd coloring}\footnote{This 
notion was introduced as \textit{odd coloring} by Petru\v{s}evski and \v{S}krekovski~\cite{PeSk2022},
but there is a previous notion of ``odd coloring'' in the literature, the difference 
being that the present one is proper and the previous one is not necessarily so.},
which is a proper coloring in which for every non-isolated vertex there is a color
appearing an odd number of times among its neighbors.
The least~$k$ such that~$G$ has a proper odd $k$-coloring is called the
\defi{proper odd chromatic number} of~$G$ and is denoted by~$\chio(G)$.

We clearly have $\chi(G) \le \chio(G) \le \chipcf(G)$, where $\chi(G)$ is the
proper chromatic number of~$G$.
We say that a graph class~$\cG$ is \defi{linearly $\chipcf$-bounded} (\defi{linearly
$\chio$-bounded}) if there is a linear function~$f$ such that 
$\chipcf(G) \leq f(\chi(G))$ ($\chio(G) \leq f(\chi(G))$) for every $G \in \cG$.

\begin{problem}[Caro \textit{et al.}~\cite{CaPeSk2023}]
\label{prob:chi}
    Find ``generic'' graph classes that are linearly $\chipcf$-bounded.
\end{problem}

As a starting point for their problem, Caro \textit{et al.}~\cite{CaPeSk2023} showed that the class
of claw-free graphs (graphs excluding $K_{1,3}$ as an induced subgraph) is linearly
$\chipcf$-bounded.  Namely, they showed that every claw-free graph~$G$ satisfies
$\chipcf(G) \le 2\Delta(G)+1$.  Since every claw-free graph satisfies $\Delta(G)/2 \le \chi(G)$,
this gives $\chipcf(G) \le 4\chi(G)+1$.  This paper improves this bound for claw-free (and more so
for quasi-line) graphs to a near-tight bound.  Moreover, we show that convex-round graphs and
permutation graphs are linearly $\chipcf$-bounded, and also give a class of circle graphs that are
$\chio$-bounded but not $\chipcf$-bounded.  We present these last results in
Section~\ref{sec:introbip}, where we also present a lemma that reduces the problem of deciding if a
hereditary class is linearly $\chipcf$-bounded to deciding if the bipartite graphs in the class are
$\chipcf$-bounded by an absolute constant.  Our results on claw-free graphs and quasi-line graphs
are independent of this lemma and also give evidence to a conjecture of Caro \textit{et
  al.}~\cite{CaPeSk2023}; these are presented in Section~\ref{sec:introclaw}.

\subsection{Boundedness and bipartite graphs}
\label{sec:introbip}

Regarding negative results concerning $\chio$-boundedness, it is known~\cite{PeSk2022}
that the \defi{full subdivision} $S(G)$ of a graph~$G$ (i.e., the graph obtained from~$G$
by subdividing each edge once) is a bipartite graph with $\chipcf(S(G)) \geq \chio(S(G)) \geq \chi(G)$, 
hence showing that the class of bipartite graphs is not $\chio$-bounded. 
It turns out, as we will prove, that bipartite graphs are \emph{the} obstacle 
for linear $\chio$- and $\chipcf$-boundedness. 
Taking a step in this direction and as a way of answering Problem~\ref{prob:chi},
Liu~\cite{Liu2024} recently proved the following result (indeed, for the choice version of $\chipcf$).

\begin{theorem}[Liu~\cite{Liu2024}]
\label{thm:liu}
    For every positive integer~$\ell$, there exists an integer~$c_\ell$ such that
    for every graph~$G$ with no odd $K_\ell$-minor, if $\chipcf(H) \le t$ for every
    induced bipartite subgraph~$H$ of~$G$, then $\chipcf(G) \leq t + c_\ell$.
\end{theorem}

Graphs excluding an odd minor include those excluding a minor, but can also be arbitrarily dense.
However, graphs with excluded odd minors have bounded chromatic number (see, for
example,~\cite{GGRSV09}). 
Theorem~\ref{thm:liu} essentially tells us that if~$\chipcf(G)$ is large, then
either~$G$ contains a large clique as an odd minor or has an induced bipartite
subgraph~$H$ with~$\chipcf(H)$ being large.
The following lemma (to be proved in Section~\ref{sec:bip}) guarantees that if~$\chipcf(G)$ is large, then either~$G$ 
has large chromatic number or contains an induced bipartite subgraph~$H$ 
with~$\chipcf(H)$ being large.
Hence, our result (qualitatively) generalizes Theorem~\ref{thm:liu}, while
in the latter the dependence on~$t$ is better.

\begin{restatable}{lemma}{reducetobip}
\label{lem:reducetobip}
    Let $t$ be an integer. If $G$ is a graph in which  every induced $(A,B)$-bipartite subgraph has a proper odd (resp. proper conflict-free) coloring~$c$ with $|c(A)|, |c(B)| \leq t$, then $\chio(G)$ (resp. $\chipcf(G)$) is bounded above by  $t^2\chi(G)-2t(t-1)$.
\end{restatable}

Lemma~\ref{lem:reducetobip} provides the following way to attack Problem~\ref{prob:chi}:
if we want to know whether a hereditary class~$\cG$ of graphs is linearly
$\chipcf$-bounded, then it is enough to check that there is an absolute constant~$k$
such that every bipartite graph $G \in \cG$ satisfies $\chipcf(G) \le k$.
With this approach, we obtain linearly $\chipcf$-boundedness for two classes.

A graph is a \defi{permutation graph} if it is the intersection graph of
segments having endpoints on two parallel lines.
Hickingbotham~\cite{Hickingbotham23} proved that every graph satisfies $\chipcf(G) \le 2\mathrm{col}_2(G)-1$,
where $\mathrm{col}_2(G)$ denotes the 2-coloring number of~$G$ (see~\cite{KY03} for a definition),
while \cite[Lemma 1(a)]{DPUY22} implies the bound $\mathrm{col}_2(G)\le 25\chi(G)$ if~$G$ is
a permutation graph.
Thus \hbox{$\chipcf(G) \le 50\chi(G)-1$} for every permutation graph~$G$.
However, we are able to improve this as follows.

\begin{restatable}{theorem}{perm}
\label{thm:perm}
    For every permutation graph~$G$, we have $\chipcf(G) \leq 3\chi(G)$.
\end{restatable}

A graph is called \defi{convex-round} if its vertices have a circular ordering such that
every neighborhood is an interval~\cite{BaHuYe2000}.
Note that the result of Hickingbotham cannot be applied to convex-round graphs,
since balanced complete bipartite graphs are convex-round and have unbounded
2-coloring number.
However, Lemma~\ref{lem:reducetobip} allows us to present a new linearly 
$\chipcf$-boundedness result for this class.

\begin{restatable}{theorem}{convex}
\label{thm:convex}
    For every (non-trivial) convex-round graph~$G$, we have $\chipcf(G) \leq 9\chi(G)-12$.
\end{restatable}

To prove this result we rely on the fact that any $(A,B)$-bipartite graph~$G$ 
that is convex-round is \defi{biconvex}~\cite{BaHuYe2000}.
This means, for $V \in \{A, B\}$ and $W \in \{A, B\} \setminus \{V\}$, there is 
a linear order~$L$ of~$V$ such that for every $w \in W$ the vertices in~$N(w)$ 
appear consecutively, as an interval $I_L(w)$, in~$L$.
Theorem~\ref{thm:convex} thus follows from Lemma~\ref{lem:reducetobip} 
and the following result (to be proved in Section~\ref{sec:bip}).

\begin{restatable}{lemma}{lembiconvex}
\label{lem:biconvex}
    Every biconvex graph~$G$ with bipartition $(A,B)$ has a proper conflict-free 
    $6$-coloring~$c$ such that $|c(A)|, |c(B)| \leq 3$.
\end{restatable}

Indeed, any bipartite graph that is a permutation graph 
is also biconvex~\cite{SpBrSt1987}.
So using Lemma~\ref{lem:reducetobip} and Lemma~\ref{lem:biconvex} we also obtain
$\chipcf(G) \leq 9\chi(G)$ for any permutation graph~$G$, a bound which is improved
by Theorem~\ref{thm:perm}.
In Proposition~\ref{prop:lb} we show that neither convex-round nor permutation 
graphs can have the identity as their $\chipcf$-bounding function.

We then look at the limits of $\chio$- and $\chipcf$-boundedness, obtaining results
that in particular tell us that Theorem~\ref{thm:perm} and Lemma~\ref{lem:biconvex} 
cannot be extended in very natural ways. 
For instance, since permutation graphs are the comparability graphs of $2$-dimensional posets,
it is natural to wonder whether Theorem~\ref{thm:perm} extends to higher-dimensional posets.
However, for every natural~$r$ there is a graph~$G$ with $\chi(G) \geq r$ and 
the dimension of the poset whose comparability graph is~$S(G)$ is 
at most~$4$~\cite{OsRo2005,TrWa2014}.
Hence, with Proposition~\ref{prop:chipcf-sg}, we have bipartite comparability graphs 
of $4$-dimensional posets of unbounded proper odd chromatic number.
However, if the dimension of the poset whose comparability graph is~$S(G)$ is at most~$3$,
then~$G$ is planar~\cite{Sc1989} and, in particular, has bounded chromatic number and, 
by a result of Ahn, Im, and Oum~\cite{AhnImOum2022}, also~$S(G)$ has bounded proper conflict-free chromatic number.
We show below that comparability graphs of $3$-dimensional posets are not $\chipcf$-bounded or even $\chio$-bounded. 

Another class containing permutation graphs is the class of \defi{circle graphs},
i.e. intersection graphs of chords of a circle. 
For their part, \defi{segment intersection graphs}, i.e. intersection graphs of horizontal 
and vertical segments in the plane, contain bipartite permutation graphs.
 On the other hand, a ``one-sided'' generalization of biconvex is the notion of \defi{convex bipartite}, i.e., 
for one of $A, B$, say~$A$, there is a linear order~$L$ such that for every $w \in B$ the vertices
in~$N(w)$ appear consecutively, as an interval $I_L(w)$, in~$L$~\cite{LP81}.
The following theorem establishes some limits for $\chipcf$- and 
$\chio$-boundedness: it tells us that Theorem~\ref{thm:perm} cannot be extended to 
$\chipcf$-boundedness of 3-dimensional posets or circle graphs, that Lemma~\ref{lem:biconvex} 
cannot be extended to convex bipartite graphs or bipartite segment intersection graphs, and 
it separates proper odd from proper conflict-free chromatic number.

\begin{theorem}\label{thm:unbounded}
    \qitem{(a)}There is a class $\cG$ of convex bipartite grid-intersection graphs that are 
    comparability graphs of $3$-dimensional posets with unbounded proper odd chromatic number.
    
    \qitem{(b)}There is a class $\mathcal{H}$ of bipartite circle graphs of bounded proper odd chromatic 
    number and unbounded proper conflict-free chromatic number.
\end{theorem}

Our construction in Theorem~\ref{thm:unbounded} (b) is identical with one from~\cite{KRS21}
(given in a different representation), where it is shown that even the (improper) conflict-free
chromatic number of circle graphs is unbounded.
However, here we emphasize that this class separates proper odd and proper conflict-free chromatic 
number.
Whether circle graphs are $\chio$-bounded remains open, see Question~\ref{quest:circle}. 

To conclude our study of bipartite graphs, we make the following contribution 
regarding proper conflict-free colorings of full subdivisions.
As mentioned earlier, the full subdivision~$S(G)$ of a graph~$G$
satisfies $\chipcf(S(G)) \geq \chio(S(G)) \geq \chi(G)$.  
On the other hand, Ahn, Im, and Oum~\cite{AhnImOum2022}, as cited above, showed that
$\chipcf(S(G)) \leq \max(5,\chi(G))$ for every graph~$G$ (as a way of proving that 
it is NP-complete to decide if a bipartite graph~$H$ has $\chipcf(H)\le k$ for $k\ge 5$).
We study the case of small~$\chi(G)$.

\begin{restatable}{proposition}{subdivisions}
    \label{prop:chipcf-sg}
    If~$G$ is a graph and~$S(G)$ is its full subdivision, we have that 
    \begin{enumerate}[(i)]
         \item\label{it:chipcf-sg-chi-3} if $\chi(G) \leq 3$, then $\chipcf(S(G)) \leq 4$,
        \item\label{it:chipcf-sg-chi-4} if $\chi(G) \leq 4$, then $\chio(S(G)) \leq 4$.
    \end{enumerate}
    There exist bipartite graphs~$G$ with $\chio(S(G)) = 4$.
\end{restatable}

Note that the above results are all tight. 
However, we do not know whether the bound $\chipcf(S(G))\leq 5$ is tight 
when $\chi(G)=4$, see Question~\ref{quest:4,5}.

In the next section we mention our results which improve the linear 
$\chipcf$-bounding function for claw-free and quasi-line graphs.
For these classes our main results are nearly tight, and cannot be obtained 
through Lemma~\ref{lem:reducetobip}, so we use more specific tools for these classes.

\subsection{Claw-free and quasi-line graphs}
\label{sec:introclaw}

The study of (improper) conflict-free edge coloring dates back to the work of
Pyber~\cite{Py1991} who showed that, for every simple graph~$G$, its line graph can be
conflict-free colored with at most~$4$ colors.
This, together with Vizing's Theorem, tells us that the line graph of every simple 
graph~$G$ can be properly conflict-free colored with $4\Delta(G) + 4$ colors.

Going beyond line graphs, (improper) conflict-free coloring has been considered for 
the class of claw-free graphs.
For instance, Bhyravarapu, Kalyanasundaram, and Mathew~\cite{BhKaMa2021} proved that
every claw-free graph~$G$ can be conflict-free colored with $O(\log (\Delta (G)))$ colors,
thus extending a result of D\c{e}bski and Przyby{\l}o~\cite{DePr2022}.
For proper conflict-free colorings, Caro \textit{et al.}~\cite{CaPeSk2023} proved that,
for every claw-free graph~$G$, we have $\chipcf(G) \le 2\Delta(G) + 1$.
For line graphs, this improves (only) the additive constant of the bound obtained
through Pyber's result and Vizing's Theorem (since if~$H$ is the line graph of 
a simple graph~$G$, then we could have $\Delta(H) = 2\Delta(G) - 2$).

Note that for every claw-free graph~$G$ we have $\chipcf(G)=\chio(G)$, because each
neighborhood has independence number at most two, and thus a color appearing 
an odd number of times in a neighborhood must appear exactly once. 
Therefore, two very recent results of Dai, Ouyang, and Pirot~\cite{DaiOP24} provide 
different improvements on the above bound of Caro \textit{et al}.
They proved that any graph~$G$ with maximum degree~$\Delta$ satisfies both 
$\chio(G) \le \Delta + O(\ln \Delta)$ 
and $\chio(G) \le \lfloor \frac{3\Delta}2 \rfloor + 2$ for every $\Delta$.

We improve all these bounds for claw-free graphs as follows. 

\begin{restatable}{theorem}{claw}~\label{thm:claw}
    Every claw-free graph~$G$ satisfies $\chipcf(G)\le \Delta(G)+6$.
\end{restatable} 

Recall that this implies $\chipcf(G) \leq 2\chi(G)+6$ for claw-free graphs.
Our proof is constructive, and a polynomial-time algorithm for finding the coloring can be derived from it.

Besides being an improvement on the upper bounds cited above, Theorem~\ref{thm:claw} also
gives support to the existence of a constant~$C$ such that $\chipcf(G) \leq \Delta(G) +C$ for each
graph~$G$, a belief expressed by Caro \textit{et al.}~\cite{CaPeSk2023}, who conjectured
that such a constant should be~1.
This conjecture has received considerable attention recently:
Caro \textit{et al.} proved that every graph~$G$ satisfies 
$\chipcf(G) \le 5\Delta(G)/2$. Cranston and Liu~\cite{CranstonL24} proved that if $\Delta(G) \ge 10^8$, then $\chipcf(G) \le \lceil 1.6550826\,\Delta(G) + \sqrt{\Delta(G)} \rceil$,
with slightly weaker bounds when $\Delta(G) \ge 750$. Most recently, Liu and Reed~\cite{LiuReed25} 
proved that the conjecture of Caro \textit{et al.} holds asymptotically.  They
proved that for~$\Delta$ large enough, every graph with maximum degree at most~$\Delta$ 
has a proper conflict-free coloring with $(1 + o(1))\Delta$ colors. 

A graph~$G$ is said to be a \defi{quasi-line graph} if, for every $v \in V(G)$, the set 
$N(v)$ can be expressed as the union of two cliques.
The class of quasi-line graphs is a proper superset of the class of line graphs
(of multigraphs) and a proper subset of the class of claw-free graphs. 
It is also a superset of the class of \defi{concave-round graphs} which are the
complements of convex-round graphs considered in Theorem~\ref{thm:convex}, see~\cite{Saf20}.
Chudnovsky and Ovetsky~\cite{ChOv2007} proved that if~$G$ is a quasi-line graph, 
then $\chi(G) \le \frac 32\omega(G)$, generalizing a classic result by Shannon~\cite{Sh1949}.
We further improve Theorem~\ref{thm:claw} for quasi-line graphs as follows.

\begin{restatable}{theorem}{quasi}~\label{thm:quasi}
    Every quasi-line graph~$G$ satisfies $\chipcf(G)\le \Delta(G)+4$.
\end{restatable}  

Since quasi-line graphs satisfy $\Delta(G) \le 2\omega(G)-2$, any such graph~$G$ satisfies $\chipcf(G) \le 2\chi(G) + 2 \enspace$.

\subsection{Organization of the paper and notation}
\label{sec:introorg}

The rest of the paper is organized as follows.
In Section~\ref{sec:claw}, we prove Theorems~\ref{thm:claw} and~\ref{thm:quasi}.
In Section~\ref{sec:subdiv}, we prove Proposition~\ref{prop:chipcf-sg}.
In Section~\ref{sec:reduce}, we prove Lemmas~\ref{lem:reducetobip} and~\ref{lem:biconvex} which are used in the proofs of Theorems~\ref{thm:perm} and~\ref{thm:convex}. 
Sections~\ref{sec:biconvex} and~\ref{sec:circle} are devoted to prove Theorem~\ref{thm:unbounded}\,(a) and Theorem~\ref{thm:unbounded}\,(b). 

To ease the readability of this paper, we adopt a convention regarding  the notation used to
refer to certain types of objects such as colorings of a graph, colors of vertices, and vertices of a graph.
  Throughout the paper, for colorings, we use lowercase letters
such as $c$, $c'$, $c''$, $c_1$, $c_2$ (letter~$c$ with different super or subscripts); for colors, we use
lowercase Greek letters such as $\alpha$, $\beta$, $\gamma$ and also positive integers in a set
$[k] = \{1, 2,\ldots, k\}$, for some $k$; for vertices, we use lower case letters such as $u$, $v$,
$w$, $y$, $z$ (possibly with subscripts). 

If~$c$ is a coloring of a graph~$G$ and $v$ is a vertex of $G$, then we denote by $W_o(c,v)$
(resp. $W_{pcf}(c,v)$) the set of colors that appear an odd number of times (resp. exactly once) in
the neighborhood of~$v$.  We denote by~$U(c)$ the set of all vertices $v$ in~$G$ with $W_o(c,v)$
empty. If $Z$ is a subset of $V(G)$, then $c(Z)$ denotes the set of colors of the coloring~$c$ that
appears in~$Z$.

For convenience, here a \emph{clique} of a graph~$G$  may refer either to a complete
subgraph of~$G$, or to a set of vertices $K$ of $V(G)$ such that $G[K]$ induces a complete subgraph of~$G$.

We write $N[u]$ for the \defi{closed neighborhood} of~$u$, defined as $N(u) \cup \{u\}$.

\section{Claw-free graphs and Quasi-line graphs}
\label{sec:claw}

In this section, we prove that any claw-free graph~$G$ has a proper conflict-free (equivalently, odd) $(\Delta(G)+6)$ coloring. 
We also improve to $\Delta(G)+4$ the previous upper bound for the class of quasi-line graphs. 
Subsection~\ref{subsec:genprelim} contains general terminology and statements regarding 
odd colorings in graphs that are used in the proofs and in the rest of the paper.

\subsection{Preliminaries}
\label{subsec:genprelim}


For each vertex~$v$ and each color $\ci \notin c(N[v])$, we denote by~$c^\ci_v$ the coloring 
obtained from~$c$ by recoloring~$v$ with color~$\ci$.
Clearly, if $c$ is a proper coloring, as $\ci \notin c(N[v])$, the coloring~$c^\ci_v$ is also a proper coloring of~$G$. 

Note that, due to parity, if~$c$ is any proper coloring of a graph
and $v$ is a vertex of odd degree, then 
$W_o(c,v)$ is non-empty but if $v$ has even degree, then $W_o(c,v)$ may be empty.
Moreover, $|W_o(c,v)|$ and~$d(v)$ have the same parity.


In the proof of our main result of this section, for a given graph~$G$, we construct a sequence 
of proper colorings $c$ such that their associated sets $U(c)$ decrease (with respect to inclusion) along the 
sequence. 
The sequence starts with an arbitrary proper $(\Delta(G)+6)$-coloring $c$ of~$G$ and continues 
unless  $U(c)$ is empty. 
To obtain the next coloring $c'$, we take any vertex~$v$ with 
$W_o(c,v)$ empty  and recolor one of its neighbors~$u$ with a color $\ci \notin c(N[u])$.
By definition, this coloring is~$c^\ci_u$.
We shall prove that
$W_o(c^\ci_u,v)$ is non-empty.
In fact, we get a bit more.

\begin{lemma}
\label{l:emptytocritical}
    Let~$G$ be a graph and~$c$ be a proper coloring of~$G$.
    For each vertex~$u \in V(G)$ and each color $\ci \notin c(N[u])$, we have that 
    $$U(c) \cap U(c^\ci_u) \cap N(u) = \emptyset\,.$$
\end{lemma}
\begin{proof}
    Let $v \in U(c) \cap N(u)$.
    Then we know that each color~$\cb$ is used by the coloring~$c$ an even number of times
    in~$N(v)$.
    Particularly, $c(u)$ and~$\ci$ are used by coloring~$c$ an even number of times in~$N(v)$. 
    The colorings~$c$ and~$c^\ci_u$ differ only at the color of vertex~$u$.
    Then, colors~$c(u)$ and~$\ci$ are used in~$N(v)$ by the coloring~$c^\ci_u$ an odd number of 
    times.
    Hence $W_o(c^\ci_u,v) = \{c(u),\ci\}$, and thus $v \notin U(c^\ci_u)$. 
\end{proof}

If the set $U(c^\ci_u) \setminus U(c)$ is empty, then the new desired coloring is~$c^\ci_u$.
The main technical obstacle that we have to overcome is dealing with vertices in 
$U(c^\ci_u) \setminus U(c)$.
The reason why a vertex~$w$ belongs to~$U(c^\ci_u)$ but not to~$U(c)$ is the following:
$w$ is a neighbor of~$u$ such that colors~$c(u)$ and~$\ci$ are the only colors appearing 
an odd number of times among its neighbors; namely, $W_o(c,w) = \{c(u),\ci\}$.
By the choice of color~$\ci$, no neighbor of~$u$ has color~$\ci$.
Hence, some neighbor~$z$ of~$w$ has color~$\ci$, and~$z$ and~$u$ are not adjacent. 
This motivates the following definition.

A vertex~$w$ is \defi{$c$-critical} if it has two neighbors, $w_1$ and~$w_2$, 
with $w_1 \notin N(w_2)$, $c(w_1) \neq c(w_2)$, and $W_o(c,w) = \{c(w_1),c(w_2)\}$.
A $c$-critical vertex creates a problem because it has only two colors in its neighborhood 
appearing an odd number of times.
However, the same fact limits the number of colors that appears in its neighborhood,
as proved in the next result. 

\begin{lemma}
\label{l:criticaldeg}
    Let $c$ be a proper coloring of a graph $G$.
    If $w$ is a $c$-critical vertex, then $|c(N[w])| \leq d(w)/2 + 2$.
\end{lemma}
\begin{proof}
    Directly from the definition of $c$-critical: $|W_o(c,w)|$ and $d(w)$ have the same parity,
    and colors in $c(N[w]) \setminus W_o(c,w)$ appear in at least two neighbors of~$w$. 
\end{proof}

A color~$\ci$ is \defi{$c$-safe} for a vertex~$u$ if $\ci \notin c(N[u])$ and 
$W_o(c,w) \neq \{c(u),\ci\}$ for each $w \in N(u)$. 

\begin{lemma}
\label{l:useofsafe}
    Let~$c$ be a proper coloring of a graph $G$.
    Let~$u$ be a vertex and $\ci \notin c(N[u])$.
    The color~$\ci$ is $c$-safe for~$u$ if and only if $U(c^\ci_u) \setminus U(c)$ is empty. 
    Moreover, $w \in U(c^\ci_u) \setminus U(c)$ if and only if $W_o(c,w) = \{c(u),\ci\}$ 
    and $w \in N(u)$. 
\end{lemma}
\begin{proof}
    Let us assume that~$\ci$ is not $c$-safe for~$u$.
    Then, there is a neighbor~$w$ of~$u$ such that $W_o(c,w) = \{c(u),\ci\}$.
    This implies that~$c(u)$ and~$\ci$ are the only colors in $c(N(w))$ that appear an odd 
    number of times.
    In coloring~$c^\ci_u$, the color~$c(u)$ is replaced by color~$\ci$, and then both colors 
    appear an even number of times in $c^\ci_v(N(w))$.
    Hence, $W_o(c^\ci_u,w)$ is empty and, therefore, $w \in U(c^\ci_u) \setminus U(c)$.

    Conversely, let $w \in U(c^\ci_u) \setminus U(c)$.
    Then $W_o(c^\ci_u,w)$ is empty, and $W_o(c,w)$ is not empty.
    Since $W_o(c^\ci_u,z) = W_o(c,z)$ for each $z \notin N(u)$, we get that $w \in N(u)$.
    Let us assume that there is a color $\cb \in W_o(c,w) \setminus \{c(u),\ci\}$.
    But then we would have~$\cb$ in $W_o(c^\ci_u,w)$, since~$c^\ci_u$ and~$c$ differ only in the
    colors~$c(u)$ and~$\ci$.
    Hence, $W_o(c,w) \subseteq \{c(u),\ci\}$.
    If some color $c'\in \{c(u),\ci\}$ does not belong to $W_o(c,w)$, then it appears an even 
    number of times in~$N(w)$.
    But as $u \in N(w)$ and it changes from color~$c(u)$ to color~$\ci$ in $c^\ci_u$, 
    the color~$\cb$ would appear an odd number of times in $c^\ci_u(N(w))$, which is not
    possible, as $W_o(c^\ci_u,w)$ is empty. 
    We conclude that $W_o(c,w) = \{c(u),\ci\}$.
\end{proof}

From Lemma~\ref{l:useofsafe}, we have that for each vertex~$u$ and each color
$\ci \notin c(N[u])$, the vertices in the set $U(c^\ci_u) \setminus U(c)$ are $c$-critical.
In fact, if $w \in U(c^\ci_u) \setminus U(c)$, then $W_o(c,w) = \{c(u),\ci\}$.
Hence, there is a neighbor~$z$ of~$w$ with $c(z) = \ci$.
Since $\ci \notin c(N[u])$, we get that~$z$ is not adjacent to~$u$.

\begin{lemma}
\label{l:safeandcritical}
    Let~$c$ be a proper coloring of a graph~$G$ and let $u \in V(G)$.
    For each $\ci,\cb \notin c(N[u])$, $\ci \neq \cb$, the sets $U(c^\ci_u) \setminus U(c)$ and 
    $U(c^\cb_u) \setminus U(c)$ are disjoint.
    Moreover, for each $\ci \neq \cb$, each vertex in $U(c^\cb_u) \setminus U(c)$ is
    $c^\ci_u$-critical.
\end{lemma}
\begin{proof}
    Let~$w_\ci, w_\cb \in N(u)$ such that $w_\ci \in U(c^\ci_u) \setminus U(c)$ and 
    $w_\cb \in U(c^\cb_u) \setminus U(c)$.
    From Lemma~\ref{l:useofsafe} we know that $W_o(c,w_\ci) = \{c(u),\ci\}$ and 
    $W_o(c,w_\cb) = \{c(u),\cb\}$.
    Since $\ci,\cb \notin c(N[u])$ and $\ci \neq \cb$, we get that $w_\ci \neq w_\cb$. 

    Now, with respect to the coloring~$c^\ci_u$, we have that $W_o(c^\ci_u,w) = \{\ci,\cb\}$
    for each neighbor~$w$ of~$u$ with $W_o(c,w) = \{c(u),\cb\}$. 
    Thus, $w$ is $c^\ci_u$-critical.
\end{proof}

From the previous lemmas, we have that if~$\ci$ is not $c$-safe for~$u$, then 
$U(c^\ci_u) \setminus U(c)$ is non-empty and $(U(c^\ci_u) \cap U(c^\cb_u)) \setminus U(c)$ is 
empty, if $\ci \neq \cb$.
Moreover, the set 
\begin{equation*}
    \bigcup_{\ci \notin c(N[u])} (U(c^\ci_u)\setminus U(c))
\end{equation*}
is a subset of~$N(u)$ and it contains only $c$-critical vertices.

\subsection{Claw-free graphs}

A graph~$G$ is \defi{claw-free} if every vertex $v \in V(G)$ does not have an independent 
set of size~3 in~$N(v)$.
Then, a graph~$G$ is claw-free if and only if for every two adjacent vertices~$v$ and~$u$, 
$N(u) \setminus N[v]$ is a clique.

Notice that in any proper coloring~$c$ of a claw-free graph, any color~$\ci$ appears 
at most twice in the neighborhood of each vertex. 
Then, for each vertex~$v$, a color $\ci \in W_o(c,v)$ if and only if there is exactly one neighbor
of~$v$ with color~$\ci$.
Therefore, a coloring (proper or not) of a claw-free graph is odd if and only if it is conflict-free.

\begin{lemma}
\label{l:criticalclawfree}
    Let~$c \colon V(G) \to [\Delta(G)+6]$ be a proper coloring of a claw-free graph~$G$.
    If a vertex~$v$ has~$t$ neighbors which are $c$-critical, then it has at least
    $\min\{2,t\}$ neighbors which are $c$-critical and have a $c$-safe color. 
\end{lemma}
\begin{proof} 
    We first prove that if a $c$-critical neighbor of~$v$ does not have a $c$-safe color,
    then it has at least four neighbors which are $c$-critical and also neighbors of~$v$. 
    This, in particular, proves the statement when $t = 1$.

    Let~$z$ be a $c$-critical neighbor of~$v$ and let $\Delta = \Delta(G)$. 
    By Lemma~\ref{l:useofsafe}, we know that for each color~$\ci$ not in~$c(N[z])$ there is a vertex 
    $z_\ci \in U(c^\ci_z) \setminus U(c) \subseteq N(z)$ with $W_o(c,z_\ci) = \{c(z),\ci\}$.
    As~$c$ is a $(\Delta+6)$-coloring and~$z$ is $c$-critical, Lemma~\ref{l:criticaldeg} 
    implies that there are at least $(\Delta+6) - d(z)/2 - 2$ colors not in~$c(N[z])$,
    which says that~$z$ has at least $\Delta + 4 - d(z)/2$ neighbors that are $c$-critical.

    We now prove that at least four of these $c$-critical neighbors of~$z$ 
    also belong to~$N(v)$.
    For this purpose, it suffices to prove that the set~$K$ given by
    \begin{equation*}
        K := N(z) \setminus N[v] 
    \end{equation*}
    contains at most $d(z)/2$ vertices, since then at least $\Delta + 4 - d(z) \geq 4$
    of the $c$-critical neighbors of~$z$ also belong to~$N(v)$. 

    We use the above-mentioned fact that, in a claw-free graph, the set~$K$ forms a clique
    in~$G$; and let $W_o(c,z) = \{c(z_1), c(z_2)\}$, with $z_1,z_2 \in N(z)$ and $z_1 \notin N(z_2)$.
    This implies, on the one hand, that at most one of~$z_1$ or~$z_2$ belongs to~$K$,
    as they are non-adjacent. 
    On the other hand, it implies that each vertex of~$K$ must have a different color in 
    the coloring~$c$.
    Then, $K$ can contain at most one vertex of those whose color appears twice in~$N(z)$, 
    and one of the vertices~$z_1$ or~$z_2$.
    Therefore, $|K| \leq |c(N[z])| - 2 = d(z)/2$.

    We now prove that when $t \geq 2$, vertex~$v$ has at least two $c$-critical neighbors,
    both of which have a $c$-safe color.

    Let~$z$ be any $c$-critical neighbor of~$v$ and let~$C$ be the set of all neighbors 
    of~$v$ which are $c$-critical, except for~$z$.
    Then, as $t \geq 2$, we get that~$C$ is not empty.
    We claim that~$C$ has a vertex that has a $c$-safe color.
    For a contradiction, let us assume that no vertex in~$C$ has a $c$-safe 
    color.
    
    Let us define the digraph $D$ with $V(D) = C$ and such that $(x,y) \in E(D)$ if and only if $c(x) \in W_o(c,y)$.
    For each $x \in C$, we know that its in-degree in $D$ is at most two since $|W_o(c,x)|=2$. 
    Hence, $| E(D)| \leq 2|C|$. 

    We have shown that if a $c$-critical neighbor of~$v$ does not have a $c$-safe color,
    then it has at least four neighbors which are $c$-critical and belong to $N(v)$.
    Hence, under our assumption that no vertex in~$C$ has a $c$-safe color, we get that
    each vertex in~$C$ has out-degree at least three in~$D$.
    This implies that $|E(D)| \geq 3|C|$, which is in contradiction with $|E(D)| \leq 2|C|$ since, 
    by hypothesis, $C$ is not empty.

    To finish the proof, we repeat the previous argument but now taking as~$z$ 
    a $c$-critical neighbor of~$v$ having a $c$-safe color whose existence we already proved.
\end{proof}

\begin{lemma}
\label{l:lowhappinesisfine}
    Let~$c \colon V(G) \to [\Delta(G)+6]$ be a proper coloring of a claw-free graph~$G$.
    Let~$v$ be a vertex such that no vertex in $N[v]$ is $c$-critical. 
    Then, each $w \in N(v)$ with $|W_o(c,w)| \leq 4$, has a $c$-safe color.
\end{lemma}
\begin{proof}
    For any vertex~$w \in V(G)$, we have that 
    \begin{equation*}
        |c(N[w])| = \frac{d(w) - |W_o(c,w)|}{2} + |W_o(c,w)| + 1 = \frac{d(w) + |W_o(c,w)|}{2} + 1\,.
    \end{equation*}
    Moreover, any clique contained in~$N(w)$ has size at most $|c(N[w])|-1$. 
    As we are assuming that~$G$ is claw-free, the set $N(w) \setminus N[v]$
    forms a clique, which implies that it has at most $|c(N[w])|-1$ vertices. 
    Furthermore, for each $\ci \notin c(N[w])$, the set $U(c^\ci_w) \setminus U(c)$ 
    is contained in $N(w) \setminus N[v]$, because each vertex 
    in $U(c^\ci_w) \setminus U(c)$ is $c$-critical and no vertex in $N[v]$ is $c$-critical,
    by hypothesis. 

    Now, if no color is $c$-safe for~$w$, then we get that the set $N(w) \setminus N[v]$
    has at least $\Delta(G) + 6 - |c(N[w])|$ vertices, one for each color not in $c(N[w])$.
    Therefore, if no color is $c$-safe for~$w$, then we get the contradiction
    \begin{equation*}
        \Delta(G)+6 \leq |N(w) \setminus N[v]| + |c(N[w])| \leq 
        2|c(N[w])| - 1 = d(w) + W_o(c,w) + 1 < \Delta(G)+6\,.
    \end{equation*}
    The result follows.
\end{proof}

Finally, we can state our main result for this section.

\claw*
\begin{proof}
    Let~$c \colon V(G) \to [\Delta(G)+6]$ be a proper coloring of~$G$.
    We shall prove that if $|U(c)| > 0$, then there is another coloring~$c'$ 
    with $|U(c')| < |U(c)|$.

    Let $v \in U(c)$.
    Suppose that there exists a neighbor~$u$ of~$v$ having a $c$-safe color~$\ci$.
    By Lemma~\ref{l:useofsafe}, $U(c^\ci_u) \setminus U(c) = \emptyset$.
    As $v \in U(c) \cap N(u)$, then by Lemma~\ref{l:emptytocritical}, $v \notin U(c^\ci_u)$. 
    Thus, $|U(c^\ci_u)| \leq |U(c)|-1$.
    Hence, we can continue under the assumption that no neighbor of~$v$ has a $c$-safe
    color.  

    Let $u \in N(v)$.
    Note that every $\ci \notin c(N[u])$ is not $c$-safe for~$u$, which means, 
    by Lemma~\ref{l:useofsafe}, that $U(c^\ci_u) \setminus U(c)$ is not empty.
    Also note that the number of colors not appearing in $c(N[u])$ is
    $\Delta(G)+6 - |c(N[u])| \geq \Delta(G)+6 - (\Delta(G)+1) = 5$.
    Fix some $\ci \notin c(N[u])$.
    For $\cb \notin c(N[u])$, with $\cb \neq \ci$, we know, by Lemma~\ref{l:safeandcritical}, 
    that all vertices in $U(c^\cb_u) \setminus U(c)$ are $c^\ci_u$-critical and that
    $U(c^\cb_u) \setminus U(c)$ is disjoint from $U(c^\ci_u) \setminus U(c)$.
    Since $\bigcup_{\cb \notin c(N[u])} (U(c^\cb_u)\setminus U(c)) \subseteq N(u)$,
    we conclude that~$u$ has at least~4 neighbors which are $c^\ci_u$-critical.

    Thus, from Lemma~\ref{l:criticalclawfree}, we obtain that~$u$ has at least two 
    $c^\ci_u$-critical neighbors both having a $c^\ci_u$-safe color.
    Then, $u$ has a $c^\ci_u$-critical neighbor~$w \neq v$ that has a $c^\ci_u$-safe 
    color~$\cg$.

    We show now that $w \notin N(v)$.
    Indeed, notice that the sets $W_o(c_u^\ci,w)$ and $W_o(c,w)$ differ exactly in colors~$c(u)$ 
    and~$\ci$, when they are not the same set. 
    Then, 
    \begin{equation*}
        |W_o(c,w)| \leq |W_o(c^\ci_u,w)| + 2 = 4\,.
    \end{equation*}

    Since no neighbor of~$v$ has a $c$-safe color, from Lemma~\ref{l:criticalclawfree} 
    we know that no neighbor of~$v$ is $c$-critical.
    Hence, $w$ cannot be a neighbor of~$v$; otherwise, as we have that $|W_o(c,w)|\leq 4$,
    it would have a $c$-safe color, by Lemma~\ref{l:lowhappinesisfine}.
   
    We prove now that the coloring $c' = (c^\ci_u)^\cg_w$ is such that~$U(c')$ is a proper 
    subset of~$U(c)$.
    For this purpose, we first show that $U(c_u^\ci) \setminus U(c)$ is contained in~$N(w)$. 

    Since no vertex in the set $N[v]$ has a $c$-safe color, due to Lemma~\ref{l:useofsafe}
    we know that $U(c^\ci_u) \setminus U(c) \subseteq (N(u) \setminus N[v])$.
    We already proved that $w \in N(u) \setminus N[v]$.

    We  use again the fact that when~$G$ is a claw-free graph, the set  
    $N(u) \setminus N[v]$ forms a clique to conclude that the set 
    $U(c^\ci_u) \setminus U(c)$ is a subset of~$N(w)$.
    This implies that
    \begin{equation*}
        U(c^\ci_u) = (U(c^\ci_u) \cap U(c)) \cup (U(c^\ci_u) \setminus U(c)) \subseteq 
            (U(c^\ci_u) \cap U(c)) \cup (U(c^\ci_u) \cap N(w)).
    \end{equation*}

    We also have that
    \begin{equation*}
        U(c') = (U(c') \cap U(c^\ci_u)) \cup (U(c') \setminus U(c^\ci_u)) = U(c') \cap U(c^\ci_u) \, ,
    \end{equation*}
    since, by Lemma~\ref{l:useofsafe}, we have that $U(c') \setminus U(c^\ci_u)$ is empty, as 
    the color~$\cg$ is $c^\ci_u$-safe for~$w$. 

    Thus, 
    \begin{equation*}
        U(c') = U(c') \cap U(c^\ci_u) \subseteq U(c') \cap U(c^\ci_u) \cap U(c)\,,
    \end{equation*}
    since, by Lemma~\ref{l:emptytocritical}, we have that $U(c') \cap U(c^\ci_u) \cap N(w)$
    is empty.
    Again by Lemma~\ref{l:emptytocritical}, we know that $v \notin U(c^\ci_u)$, since 
    $\ci \notin c(N[u])$ and $v \in U(c) \cap N(u)$. 
    Therefore, $U(c')$ is a subset of $U(c^\ci_u) \cap U(c)$ which, as 
    $v \in U(c) \setminus U(c^\ci_u)$, is a proper subset of $U(c)$.
\end{proof}

\subsection{Quasi-line graphs}

A graph~$G$ is a \defi{quasi-line} graph if, for each vertex $v \in V(G)$, there is 
a partition $\{K^1_v, K^2_v\}$ of $N(v)$ into cliques.
Then, any quasi-line graph is a claw-free graph, and any line graph is a quasi-line graph.

\begin{lemma}
\label{l:critical}
    Let~$G$ be a quasi-line graph and let~$c \colon V(G) \to \mathbb{N}$ be 
    a proper coloring of~$G$.
    Let~$u$ be a $c$-critical vertex.
    Then~$d(u)$ is even and $|c(N[u])| = d(u)/2 + 2$. 
    Moreover, for each $c$-critical vertex $w \in N(u)$, we have $d(w) = d(u)$.
\end{lemma}
\begin{proof}
    To prove the first statement, remember that by definition of~$u$ being $c$-critical,  we 
    have that $W_o(c,u) = \{c(u_1),c(u_2)\}$, where~$u_1$ and~$u_2$ are neighbors of~$u$ 
    which are not adjacent.
    As~$G$ is a quasi-line graph, there is a partition $\{K^1_u,K^2_u\}$ of $N(u)$ into 
    cliques and, since $u_1, u_2$ are not adjacent, we can assume that 
    $u_1 \in K^1_u$ and $u_2 \in K^2_u$. 

    As $K^1_u$ and $K^2_u$ are cliques, we know that 
    $|c(K^1_u)| = |K^1_u|$ and $|c(K^2_u)| = |K^2_u|$.
    Due to $W_o(c,u) = \{c(u_1),c(u_2)\}$, we actually have that 
    $c(K^1_u) \setminus \{c(u_1)\} = c(K^2_u) \setminus \{c(u_2)\}$ and, thus,  
    $|K^1_u| = |K^2_u|$.
    Therefore, $d(u) = 2|K^1_u|$ and $|c(N(u))| = |K^1_u| + 1 = d(u)/2 + 1$. 
    Hence, $|c(N[u])| = d(u)/2 + 2$.

    Let $w \in N(u)$ be a $c$-critical vertex.
    We just proved that the degree of any $c$-critical vertex is even 
    and that $|c(N[w])| = d(w)/2 + 2$. 
    Let~$w_1$ and~$w_2$ be non-adjacent neighbors of~$w$ such that 
    $W_o(c,w) = \{c(w_1),c(w_2)\}$.

    We can assume that the vertex~$w$ belongs to~$K^1_u$.
    Then, $K := (K^1_u \cup \{u\}) \setminus \{w\}$ is contained in~$N(w)$.
    Since~$K$ is a clique, we have that $|c(K)| = |K|$ and one of the vertices~$w_1$
    or~$w_2$ does not belong to~$K$.
    Hence, $|c(K)| \leq |c(N(w))| - 1 = d(w)/2$.
    Since $d(u)/2 = |K^1_u| = |K|$, we get that $d(u)/2 \leq d(w)/2$. 

    By interchanging the role of~$u$ and~$w$ in the previous argument, we get that 
    $d(u) \geq d(w)$, then proving the second statement.  
\end{proof}

\begin{lemma}
\label{l:safedigraph}
    Let~$G$ be a quasi-line graph and let~$c \colon V(G) \to [\Delta(G)+4]$ be a proper coloring of~$G$.
    Let~$u$ be a $c$-critical vertex.
    If~$u$ does not have a $c$-safe color, then for each partition $\{K^1,K^2\}$ 
    of~$N(u)$ into cliques and each $j \in \{1,2\}$, the vertex~$u$ has a neighbor in~$K^j$ 
    which is $c$-critical and has a $c$-safe color. 
\end{lemma}
\begin{proof} 
    Let $i \in \{1,2\}$ and let~$C$ be the subset of $c$-critical vertices in 
    $K^i \cup \{u\}$. 
    We consider the digraph $D$ with $V(D) = C$ and such that $(v,w) \in E(D)$ whenever
    $c(v) \in W_o(c,w)$. 

    If the arcs $(v,w)$ and $(z,w)$ are in~$E(D)$, then $c(v),c(z) \in W_o(c,w)$.
    As~$C$ is a clique in~$G$, we have that $c(v) \neq c(z)$ and thus we get that 
    $W_o(c,w) = \{c(v),c(z)\}$.
    But then, $w$ would not be $c$-critical, as~$v$ and~$z$ are adjacent. 
    We conclude that the in-degree of each vertex in~$D$ is at most one
    and thus $|E(D)| \leq |C|$. 

    Now, we show that every vertex in~$C$ not having a $c$-safe color has out-degree at 
    least two in~$D$. 
    Let $z \in C$ be a vertex with no $c$-safe color. 
    By Lemma~\ref{l:useofsafe}, we know that for each color~$\ci$ not in~$c(N[z])$ there is 
    a vertex $z_\ci \in C(c^\ci_z) \setminus C(c) \subseteq N(z)$.
    As~$c$ is a $(\Delta(G)+4)$-coloring and~$z$ is $c$-critical, Lemma~\ref{l:critical} 
    implies that there are $(\Delta(G)+4) - d(z)/2 - 2$ colors not in $c(N[z])$, 
    which says that~$z$ has $\Delta(G) + 2 - d(z)/2$ neighbors which are $c$-critical.  

    From Lemma~\ref{l:critical} we get that $d(z) = d(u)$ and hence, as~$z$ has $|K^j|$ 
    neighbors in $K^j \cup \{u\}$ and $|K^j| = d(u)/2 = d(z)/2$, we know that~$z$ has
    $d(z)/2$ neighbors not in $K^j \cup \{u\}$. 
    Hence, at most $d(z)/2$ of the $\Delta(G)+2 - d(z)/2$ neighbors of~$z$ which are 
    $c$-critical do not belong to~$C$. 
    Whence, at least $\Delta(G)+4 - d(z) - 2 \geq 2$ of these neighbors belong to~$C$. 
    Therefore, the out-degree of a vertex~$z$ in~$C$ not having a $c$-safe color is 
    at least two. 
    As $|E(D)| \leq |C|$, we get that at least one vertex in~$C$ has a $c$-safe color.
\end{proof}

\quasi*
\begin{proof}
    Let~$c \colon V(G) \to [\Delta(G)+4]$ be a proper coloring of~$G$.
    We shall prove that if $|U(c)| > 0$, then there is another coloring~$c'$ 
    with $|U(c')| < |U(c)|$.

    Let $v \in U(c)$.
    Suppose that there exists a neighbor~$u$ of~$v$ having a $c$-safe color~$\ci$.
    By Lemma~\ref{l:useofsafe}, $U(c^\ci_u) \setminus U(c) = \emptyset$.
    As $v \in U(c) \cap N(u)$, then by Lemma~\ref{l:emptytocritical}, $v \notin U(c^\ci_u)$. 
    Thus, $|U(c^\ci_u)| \leq |U(c)|-1$.
    Hence, we can continue under the assumption that no neighbor of~$v$ has a $c$-safe 
    color. 
 In particular, from Lemma~\ref{l:safedigraph}, we also can assume that no neighbor $u$
    of~$v$ is $c$-critical, otherwise $u$ and $v$ have a common neighbor with a $c$-safe color. 

    Let $u \in N(v)$ and $\ci \notin c(N[u])$. 
    We prove that~$u$ has at least~3 neighbors that are $c_u^\ci$-critical.
    In fact, $v$ is one of them.
    Moreover, we know that $\Delta(G)+4 - |c(N[u])| \geq 3$ and that for each
    $\cb \notin c(N[u])$, $\cb \neq \ci$, the set $U(c_u^\cb)\setminus U(c)$ is not empty and 
    by Lemma~\ref{l:safeandcritical}, that its vertices are $c_u^\ci$-critical.

    Let $\{K^1,K^2\}$ be a partition of $N(u)$ into cliques.
    We can assume that $v \in K^2$. 

Notice that $u$ is also $c^\ci_u$-critical as the colors of its neighbors do not change. By the same reason and the fact that $W_o(c^\ci_u,v)=\{c(u),\ci\}$, it has no $c^\ci_u$-safe colors, as it has no $c$-safe colors.
Hence, by Lemma~\ref{l:safedigraph} we obtain that~$u$ has a $c^\ci_u$-critical 
    neighbor~$w$ in $K^1$ having a $c^\ci_u$-safe color~$\cg$.  

    Since~$v$ does not have $c$-critical neighbors, and each vertex in 
    $U(c_u^\ci)\setminus U(c)$ is $c$-critical, we have that
    $(U(c_u^\ci)\setminus U(c)) \cap N(v)$ is empty, and thus
    $(U(c_u^\ci)\setminus U(c)) \subseteq K^1$.

    We prove now that the coloring $c' = (c^\ci_u)_w^\cg$ is such that~$U(c')$ is a proper 
    subset of~$U(c)$.
    It is clear that $U(c_u^\ci)\setminus U(c)$ is contained in $N(w)$, since both 
    $(U(c_u^\ci) \setminus U(c))$ and~$w$ are contained in the clique~$K^1$. 
    This implies that 
    \begin{equation*}
        U(c^\ci_u) = (U(c^\ci_u) \cap U(c)) \cup (U(c^\ci_u) \setminus U(c)) 
        \subseteq (U(c^\ci_u) \cap U(c)) \cup (U(c^\ci_u) \cap N(w)) \,.
    \end{equation*}
    We also have that
    \begin{equation*}
        U(c') = (U(c') \cap U(c^\ci_u)) \cup (U(c') \setminus U(c^\ci_u)) =
        U(c') \cap U(c^\ci_u) \,,
    \end{equation*}
    since by Lemma~\ref{l:useofsafe}, we have that $U(c') \setminus U(c^\ci_u)$ is empty,
    as the color~$\cg$ is $c^\ci_u$-safe for~$w$. 
    Thus, 
    \begin{equation*}
        U(c') = U(c') \cap U(c^\ci_u) \subseteq U(c') \cap U(c^\ci_u) \cap U(c) \,,
    \end{equation*}
    since by Lemma~\ref{l:emptytocritical} we have that $U(c') \cap U(c^\ci_u) \cap N(w)$ is 
    empty.
    Therefore, $U(c')$ is a subset of $U(c^\ci_u) \cap U(c)$ which, as
    $v \in U(c) \setminus U(c^\ci_u)$, is a proper subset of~$U(c)$.
\end{proof}

\section{Boundedness and bipartite graphs}
\label{sec:bip}

As mentioned earlier, one has $\chi(G) \leq \chio(G) \leq \chipcf(G)$ for every graph~$G$.
In this section, we investigate whether, for a family of graphs~$\cG$, there exists 
a function~$f$ such that $\chio(G) \leq f(\chi(G))$ or even $\chipcf(G) \leq f(\chi(G))$ 
for every graph $G \in \cG$.
Indeed, concerning proper conflict-free colorings~\cite[Question 6.2]{CaPeSk2023},
it has been asked for generic $\chipcf$-bounded families. 

\subsection{Full subdivisions and chromatic number}
\label{sec:subdiv}

First, we present a result that, in particular, shows the known result that no such
function~$f$ exists for large classes of bipartite graphs.
For this, given a graph~$G$, let~$S(G)$ be the (bipartite) graph obtained from~$G$ by subdividing every
edge of~$G$ exactly once. The graphs $S(G)$ are the first examples of graphs that are not $\chi_o$-bounded, see e.g.~\cite[Observation 2.7]{CaPeSk2023}. In what follows, we prove some more specific relations between $\chi(G)$ and $\chipcf(S(G))$ (or $\chio(S(G))$) for graphs of low chromatic number. Together with \cite[Lemma 3.3]{AhnImOum2022} these results show, in particular, that  
$\chipcf(S(G))$ is at most $\chi(G) + 2$ for every  graph~$G$.

\subdivisions*

\begin{proof}
  Given $G$, let ${S_V} = \{s_e \colon e \in E(G)\}$. We shall consider that $S_V$ is the set of
  vertices of $S(G)$ yielded by the subdivision of the edges of~$G$, and the other vertices of
  $S(G)$ are those originally in $G$. More precisely, $S(G)$ is a $(V(G), S_V)$-bipartite graph,  in
  which each vertex $s_e$ of $S_V$ has degree two and is adjacent to vertices $u$, $v$ in $V(G)$
  such that $e=uv\in E(G)$.

    We prove~\eqref{it:chipcf-sg-chi-3} by showing a proper conflict-free coloring
    $c' \colon V(G) \cup S_V \to [4]$ of~$S(G)$ from a $c \colon V(G) \to [\chi(G)]$ 
    proper coloring of~$G$.
    Recall that $\chi(G) \leq 3$.

    Let~$M$ be a maximum matching in~$G$ and denote by~$V(M)$ the set of vertices
    covered by~$M$ in~$G$.
    We start by putting $c'(u) = c(u)$ for all~$u \in V(M)$.
    Now we put $c'(s_{uv}) = 4$ for all $uv \in M$.
    Then, we put in $c'(s_{uv})$ a color from $[3] \setminus \{c(u), c(v)\}$
    for all $uv \in E(G) \setminus M$ such that $u \in V(M)$ or $v \in V(M)$. 
    If~$M$ is a perfect matching, then the coloring defined so far is a proper 
    conflict-free coloring of~$S(G)$, with color~$4$ 
    in $W_{pcf}(c',v)$, for every
    vertex $v$ of~$V(G)$. Note that both $c(u)$ and $c(v)$ are in $W_{pcf}(c',s_{uv})$.

    If~$M$ is not a perfect matching, then let $X = V(G) \setminus V(M)$ 
    be the set of vertices not covered by~$M$. Notice that $X$ is an independent set in $G$. 
    Put $c'(v) = 4$ for every $v \in X$.
    Then, for each $v \in X$, choose one $u \in N(v)$ and let $c'(s_{vu}) = c(v)$.
    For all other $x \in N(v) \setminus \{u\}$, we put in $c'(s_{vx})$ a color 
    from $[3] \setminus \{c(v), c(x)\}$.
    This guarantees that the color~$c(v)$ appears only once in the neighborhood of~$v$
    in the coloring~$c'$. Hence $c'$ is a proper conflict-free coloring. 
    This concludes the proof that $\chipcf(S(G)) \leq 4$. 
    
    We now prove~\eqref{it:chipcf-sg-chi-4}.
    We can assume that~$G$ is a connected graph with $\chi(G)=4$.
    We start by showing the following result.

 \begin{claim}
        Let~$T$ be a spanning tree of~$G$.
        Let~$c \colon V(G) \to [4]$ be a proper coloring of~$G$ and let~$c'$ be a coloring of 
        $\{s_e \colon e \notin E(T)\}$ such that $c'(s_{uv}) \in [4] \setminus \{c(u),c(v)\}$.
        Then, for each vertex~$r \in V(G)$, there is an extension of~$c$ and~$c'$
        to a proper $4$-coloring~$\hat{c}$ of~$S(G)$ such that for each 
        $w \in V(S(G)) \setminus \{r\}$ there is a color~$a_w$ 
    in $W_o(\hat{c},w)$,
        \end{claim}
    \begin{claimproof}
        Start by making $\hat{c}(u) = c(u)$ for all $u \in V(G)$ and
        $\hat{c}(s_e) = c'(s_e)$ for all $e \notin E(T)$.
        At this moment, note that $W_o(\hat{c},s_e)$ is non-empty, for all $e \notin E(T)$.
        Now let $r \in V(G)$.
        We assign to each vertex~$s_{uv}$ with $uv \in E(T)$ a color 
        in $[4] \setminus\{c(u),c(v)\}$ according to the following iterative process.

        Pick $uv \in E(T)$ with~$v \neq r$ being a leaf of~$T$.
        Let $\{a,b\} = [4] \setminus \{c(u),c(v)\}$.
        If $a \notin W_o(\hat{c},v)$, then let $\hat{c}(s_{uv}) = a$; otherwise, let
        $\hat{c}(s_{uv}) = b$.
        Update~$T$ with $T - uv$ and repeat the process.
        
        Because each vertex~$s_{uv}$ has received a color in $[4] \setminus \{c(u),c(v)\}$,
        the resulting coloring~$\hat{c}$ is a proper coloring of~$S(G)$ which extends the
        colorings~$c$ and~$c'$.

        Furthermore, note that after attributing a color to~$s_{uv}$ in the process, all 
        neighbors of~$v$ in~$S(G)$ are colored, and color~$a$ 
        belongs to $W_o(\hat{c},v)$.
        Also, since the degree of each vertex~$s_{uv}$ is two and $c(u) \neq c(v)$, both 
        colors 
        belong to $W_o(\hat{c},s_{uv})$.
        Thus, for all $w \in V(S(G)) \setminus \{r\}$ we have that $W_o(\hat{c},w) \neq \emptyset$.
    \end{claimproof}

    Let~$c \colon V(G) \to [4]$ be a proper $4$-coloring of~$G$ such that there is 
    a vertex~$r$ whose neighborhood is not monochromatic.
    Note that, because $\chi(G) = 4$, such $4$-coloring always exists.
    Now the result follows by induction on the number of blocks of~$G$.

    First, assume that~$G$ is 2-connected and take any spanning tree~$T$ of~$G$ 
    such that~$r$ is a leaf of~$T$ (just take any spanning tree of $G-r$ and add~$r$
    as a leaf).
    We define a coloring~$c'$ of the set $\{s_e \colon e \notin E(T)\}$ as follows.
    Let~$r'$ be the vertex adjacent to~$r$ in~$T$.
    As we are assuming that the neighborhood of~$r$ is not monochromatic in~$c$,
    there is a neighbor~$r''$ of~$r$ such that $c(r') \neq c(r'')$.
    We first color $s_{rr''}$ with color~$c(r')$ and each vertex~$s_{ru}$, with $u \neq r',r''$, 
    with some color in $[4] \setminus\{c(r),c(r'),c(u)\}$.
    Notice that with this coloring, the color $c(r') = c'(s_{rr''})$ appears exactly once
    among the neighbors of~$r$ in~$S(G)$ which are not~$s_{rr'}$.
    We then color each uncolored vertex~$s_{uv}$, with $uv \notin E(T)$ with an arbitrary
    color in $[4] \setminus\{c(u),c(v)\}$.

    From our previous claim, there is a proper coloring $\hat{c}$ of~$S(G)$ which extends~$c$ 
    and~$c'$ and such that 
    for each vertex $u \in V(S(G)) \setminus \{r\}$, $W_o(\hat{c},u)$ is non-empty.
   Moreover, $\hat{c}(s_{rr'}) \in [4] \setminus \{c(r),c(r')\}$.
    Since the color~$c(r')$ appears exactly once among the neighbors of~$r$ in~$S(G)$
    which are not $s_{rr'}$, we conclude that~$c(r') \in W_o(\hat{c},r)$.
    Therefore, $\hat{c}$ is a proper odd $4$-coloring of~$S(G)$.

    Now, consider that~$G$ is not 2-connected.
    Note that if~$r$ is not a cut-vertex, then we can also find a spanning 
    tree~$T$ of~$G$ where~$r$ is a leaf of~$T$ and use the same arguments as above
    to find an odd $4$-coloring of~$S(G)$.
    Thus the only vertices that have non-monochromatic neighborhoods are cut-vertices.

    Let~$B \subseteq G$ be a block of~$G$ containing exactly one cut-vertex~$r$ of~$G$.
    By the previous argument, all vertices in $N(r) \cap V(B)$ have monochromatic 
    neighborhoods (particularly, they only have color $c(r)$ in their neighborhood).
    Let $H := G-(B-r)$.
    Clearly, $\chi(H) \leq 4$ so, by the hypothesis, we know that there is a proper odd 
    $4$-coloring $\hat{c}_1$ of~$S(H)$.
    We now show how to extend~$\hat{c}_1$ to a $4$-coloring~$\hat{c}$ of~$S(G)$.
    For that, assume w.l.o.g that $\hat{c}_1(r) = c(r)$.

    Let~$\ci \in W_o(\hat{c}_1,r)$ and let us first define a coloring $c_2 \colon V(B) \to [4]$.
    Put $c_2(u) = \ci$ for every $u \in N(r) \cap V(B)$ and $c_2(u) = c(u)$ 
    for every $u \in V(B) \setminus N(r)$.
    Note that~$c_2$ is a proper $4$-coloring of~$B$, since $\ci \neq \hat{c}_1(r)$
    and each $u \in N(r) \cap V(B)$ only has color~$\hat{c}_1(r)$ in its neighborhood.
    By taking a spanning tree~$T$ of~$B$, any $4$-coloring~$c'_2$ of 
    $\{s_e \colon e \in E(B) \setminus E(T)\}$ with $c'_2(s_{uv}) \in [4] \setminus\{c_2(u),c_2(v)\}$
    and applying our previous claim, we get that there is a proper $4$-coloring~$\hat{c}_2$ 
    of~$S(B)$ such that 
    for each vertex $w \in V(S(B)) \setminus \{r\}$, the set $W_o(\hat{c}_2,w)$ is non-empty.

    At last, let~$\hat{c}$ be the proper $4$-coloring of~$S(G)$ defined 
    by~$\hat{c}_1$ and~$\hat{c}_2$.
    Since all vertices~$s_{ru}$ with $ru \in E(B)$ have color $c'_2(s_{ru}) \neq c_2(u) = \ci$,
    the color~$\ci$ 
    still belongs to $W_o(\hat{c},r)$. 
    
    Therefore, $\hat{c}$ is a proper odd  $4$-coloring of~$S(G)$.

    Finally, let us show that the upper bounds given
    in~\eqref{it:chipcf-sg-chi-3} and~\eqref{it:chipcf-sg-chi-4} are
    tight.  For that, consider the complete bipartite graph
    $G = K_{2k,2k}$, and let $V(G) = X \cup Y$, where $X$ and $Y$ are
    independent.  We show that $\chio(S(G)) \geq 4$ for all
    $k \geq 1$.  Let $c' \colon V(G) \cup S_V \to \{1,2,\ldots, p\}$
    be a proper odd coloring of $S(G)$.  Since every vertex $s_{xy}$ in
    $S_V$ has exactly two neighbors, it follows that the sets $X$ and
    $Y$ use together at least two different colors.  Suppose only one
    color, say~$1$, is used in~$X$, and a vertex $y$ in $Y$ gets
    color~$2$. Since $y$ has even degree, all subdivision vertices
    which are neighbors of $y$ cannot receive the same color. Thus two
    other colors different from~$1$ and~$2$ are needed to color these
    neighbors of~$y$, and therefore, $p\ge 4$. Now suppose, without
    loss of generality, that $X$ uses at least $2$~different colors.
    Since $G$ is a complete bipartite graph, the colors used (by $c'$)
    in~$X$ have to be different from the colors used in~$Y$. Thus,
    $\chio(S(G)) \geq 4$, and therefore, $\chipcf(S(G)) \geq 4$.  
  \end{proof}

The following remains open.

\begin{question}
\label{quest:4,5}
    Is there a graph $G$ with $\chi(G)=4$ and $\chipcf(S(G))=5$?    
\end{question}

\subsection{Reducing linear boundedness to bipartite graphs}
\label{sec:reduce}

The following result justifies that for questions on $\chio$-boundedness (as well as
for $\chipcf$-boundedness), we can restrict our attention to bipartite graphs.

\reducetobip*

\begin{proof} 
We give the proof for proper odd colorings since the proof for
proper conflict-free colorings is identical. To simplify notation, let
$\chi = \chi(G)$. Let $G$ be as in the statement of the lemma,
$c_0\colon V(G) \to [\chi]$ be a proper $\chi$-coloring of~$G$ and
assume without loss of generality that~$G$ has no isolated vertices.
Order the vertices of~$G$ by increasing color and compute a greedy
$\chi$-coloring~$c_1$ of~$G$ along this ordering.
It is not hard to see that~$c_1$ is indeed a $\chi$-coloring. Suppose $c_1$
defines color classes $V_1, V_2, \cdots, V_\chi$, where
$V_\ci = \{v \in V(G) \colon c_1(v) = \ci\}$.

Using now  the color classes $V_\ci$, let us define $\chi-1$ induced
bipartite subgraphs of $G$, on each of which we shall apply the
hypothesis of the lemma, and construct a proper coloring $c_2$ of $G$
(which may not be an odd coloring of $G$). To this end, for every
$\ci \in [\chi] \setminus \{1\}$, consider the induced bipartite
subgraph $H_\ci := G[V_{\ci-1} \cup V_\ci]$. Note that, in the 
bipartite subgraph $H_{\ci-1}$ the set $V_{\ci -1}$ may contain
isolated vertices, but the set $V_{\ci}$ does not contain isolated
vertices (a consequence of the greedy coloring).
    
Note also that $V(H_{\ci-1}) \cap V(H_{\ci}) = V_{\ci}$ (that is, each
pair of consecutive subgraphs $H_{\ci-1}$ and $H_{\ci}$ overlaps).
Now we shall define a coloring $c_2$ for the vertices in
$V_1, V_2, \ldots, V_\chi$, by considering for each
$\ci \in [\chi] \setminus \{1\}$, the bipartite subgraph $H_{\ci-1}$
  and a proper odd coloring of it, guaranteed by the hypothesis of the
  lemma. However, whenever we consider a bipartite subgraph
  $H_{\ci-1}$, we define the coloring $c_2$ only on the vertices of
    $V_{\ci-1}$ (the side of the bipartition that has the smallest
    color).  More precisely, we consider the proper odd coloring of
    each $H_{\ci-1}$ (guarantee by the hypothesis of the lemma), and
    consider its restriction to the vertex set $V_{\ci-1}$, and call
    it~$d_{\ci-1}$. Thus, $d_{\ci-1}:V_{\ci-1}\to [t]$ is a coloring
    of $V_{\ci-1}$ that has the property that \emph{for every vertex
      in $V_\ci$ there is a color appearing an odd number of times
      among its neighbors in $V_{\ci-1}$}.  (It helps if we think that 
      the edges of $H_{\ci-1}$ are directed from $V_{\ci}$ to
$V_{\ci-1}$, indicating that the property that we just stated holds.)  
Then, we use $d_{\ci-1}$ and define the coloring $c_2$
      on the vertices $v$ of $V_{\ci-1}$ as follows:
      $c_2(v) = \bigl(d_{\ci-1}(v),\ci-1\bigr)$. For vertices
      $v\in V_\chi$, we define $c_2(v)=\chi$.  Thus, $c_2$ is  defined
      for all vertices of $G$ in the following way.
$$c_2(v)=\begin{cases}
     \bigl(d_{\ci}(v),\ci\bigr) & \text{if } v \in V_{\ci} \text{ and } \ci\in  \{1, 2, \ldots, \chi-1\},\\
         \chi & \text{if } v\in V_{\chi}. \\
          \end{cases}
  $$

 Note that, for $v \notin V_{\chi}$, the second component $\alpha$  of $c_2(v)$ indicates that the  vertex $v$ is in the set $V_{\ci}$. Note also that $c_2$ attributes a color to the isolated vertices in each of the $\chi -1$ bipartite subgraphs that we considered above  (because in the hypothesis the induced bipartite subgraphs are not required to be free of isolated vertices).

For a coloring $c'$ of $G$, we recall that $W_o(c',v)$ denotes the
  set of colors that appear an odd number of times in the neighborhood
  of $v$. When $W_o(c',v)\ne \emptyset$, let us say that \emph{$c'$ is
    good for $v$}. Thus, so far we have that $c_2$ is good for
  $v\in \{V_2, \ldots, V_\chi\}$, but we do not know whether $c_2$ is
  good for $v\in V_1$.  Our aim is to define a proper coloring $c_4$ of $G$
  that is good for all vertices of $V(G)$. For that, we define first
  a  coloring $c_3$, then we obtain $c_4$ by combining  $c_2$ with $c_3$.

    Since $G$ has no isolated vertices, every $v\in V_1$ has a neighbor in some 
    other set~$V_\ci$;  associate~$v$ with the smallest such $\ci$ and denote it by $f(v)$.
    Now, for every $\ci\in [\chi] \setminus \{1\}$ consider the induced bipartite subgraph
    ${H'_\ci} := G[\bigcup_{v\in f^{-1}(\ci)} (\{v\}\cup (N(v)\cap V_\ci))]$, which by definition is free of isolated vertices.
    By hypothesis, there is a coloring 
    $d'_\ci \colon \bigcup_{v\in f^{-1}(\ci)} (N(v)\cap V_\ci) \to [t]$ that is good  for every vertex~$v\in f^{-1}(\ci)$. 
    Then, we define the coloring $c_3$ for
    all vertices of $G$, in the following way. If $v\in V_1$,  then let
    $c_3(v)=1$. If  $w \in V_\ci$, for $2\leq \ci \leq \chi$, we define
  $$c_3(w)=\begin{cases}
    d'_{\ci}(w) & \text{if } \exists_{v\in f^{^-1}(\ci)}:w\in N(v)\cap V_{\ci},\\
    1 & \text{otherwise.}\\
            \end{cases}
   $$
 Now we define the final coloring $c_4$ of $G$ as follows:
  \[c_4(v) = \bigl(c_2(v), c_3(v)\bigr)  \text{ for each }  v\in V(G).\]

 (Similarly to the case of $c_2$, it helps to  think that
  w.r.t. $c_3$ the edges of ${H'_\ci}$ are directed from the class
  that is a subset of $V_1$ to the class that is a subset of
  $V_\ci$. This indicates that for each $v\in V_1$, a color in
  $W_o(v,c_3)$, can be found among its neighbors in $V_\ci$, pointed
  by the arrows leaving $v$.)
 

  Let us prove now that $c_4$ is a proper odd coloring of $G$.  For
    that, let us prove first that $c_4$ is good for
    $v\in\{V_2, V_3, \ldots, V_\chi\}$.  Fix a color
    $\ci\in \{2,3,\ldots\chi\}$.  Since $c_4(v)=(c_2(v),c_3(v))$, for
    each vertex $v$ in $V_\ci$, we can see that $c_4(V_\ci)$ (that is,
    the coloring defined by $c_4$ on $V_\ci$) is a refinement of
    $c_2(V_\ci)$ induced by $c_3(V_\ci)$. Since $c_2$ is good for each
    vertex in $V_\ci$, a little argument using the fact that whenever
    an odd set is partitioned into any number of parts, at least one
    part is odd, leads us to the conclusion that $c_4$ is also good
    for the vertices in $V_\ci$.

    It remains to show that $c_4$ is good for the vertices in $V_1$.
    For each vertex $v \in V_1$, there is a color~$\ci > 1$ such that 
    $\ci = f(v)$. Let $V'_1(\ci) = \{v'\in V_1: v'\in f^{-1}(\ci)\}$,
    and let $V'_\ci$ be the set of neighbors of $V'_1(\ci)$ that belong to $V_\ci$.
     Then, by the definition of $c_3$ we have that
    $c_3$ is good for all vertices in $V'_1(\ci)$. Since
    $c_4(v)=(c_2(v),c_3(v))$, then $c_4(V'_\ci)$ can be seen as a
    refinement of $c_3(V'_\ci)$ induced by $c_2(V'_\ci)$. Thus, as in the previous case, we conclude that
    $c_4$ is good for the vertices in $V'_1(\ci)$.  Hence, $c_4$ is
    good for all vertices in $V_1$, and therefore 
    $c_4$ is a proper odd coloring of $G$. 

    Since $c_2$ uses at most $t$ different colors on each $V_\ci$
    for $1\leq \ci\leq \chi-1$, and uses~$1$ color on $V_{\chi}$, we
    have that $c_2$ uses at most $t(\chi-1)+1$ colors.  The
    coloring~$c_3$ uses in total at most $t$ colors on
    $\bigcup_{2\leq \ci\leq \chi}V_\ci$ and $1$ color on $V_1$. Thus,
    $c_4$ uses at most $t$ colors on $V_1$, at most $t^2$ colors on
    each $V_\ci$ for $2\leq \ci\leq \chi-1$, and at most $t$ colors on
    $V_{\chi}$.  Hence, $c_4$ uses at most $t+t^2(\chi-2)+t$ colors. 
    \end{proof}

\lembiconvex*
\begin{proof}
  Let $L_1$ and $L_2$ be orderings of $A$ and $B$, respectively, that
  exist because~$G$ is biconvex.  It suffices to show that~$A$ admits
  a $3$-coloring $c$ such that the set $W_{pcf}(c,v)$ is non-empty,
  for every vertex $v$ in~$B$.

    For every vertex $w \in B$ let 
    $I_w$ be the interval of $L_1$ of neighbors of $w$, and let
    $\cI = \{I_w \colon w \in B\}$.
    Also, for any $I \in \cI$, let $\cI_{|I} = \{I \cap I_w \colon w \in B, I \cap I_w \neq \emptyset\}$.
    We need the following claim.

    \begin{claim}
    \label{claim:3inclusionmin}
        For each interval $I \in \cI$, the set $\cI_{|I}$
        has at most two inclusion-minimal intervals.
    \end{claim}
    \begin{claimproof}
        Let $I = I_{w}$, and suppose that there are three inclusion-minimal intervals 
        $I \cap I_{w_1}, I \cap I_{w_2}, I \cap I_{w_3}$ in $\mathcal{I}_{|I}$. 
        Clearly, none of them corresponds to~$I$.
        For $i \in \{1,2,3\}$, let $v_i$ be the left endpoint of $I_{w_i}$.
        Since the intervals $I \cap I_{w_i}$ are inclusion-minimal, the~$v_i$ are mutually
        distinct, but each~$v_i$ also appears in~$I$ which imply that for each $i$, $N(v_i) \cap \{w_1,w_2,w_3,w\} = \{w_i,w\}$. 
        By the interval property in~$L_2$, it follows that~$w$ needs to be in the same interval with 
        each of the three~$w_1$, $w_2$, and $w_3$ in~$L_2$, which is not possible.
    \end{claimproof}
 
    We are ready to describe the coloring of~$A$. 
    Consider any inclusion-maximal set~$\mathcal{I}_M$ of mutually disjoint intervals 
    in~$\mathcal{I}$.
    (To obtain this, we can simply go from left to right through~$L_1$, and every time we see 
    a vertex uncovered by our set, we add an interval having this vertex as its left endpoint.)
    For each $I \in \mathcal{I}_M$, take the leftmost and rightmost vertices contained
    in a minimal interval of~$\mathcal{I}_{|I}$. 
    Collect these vertices in a set $S \subseteq A$ and color them alternating from left
    to right with colors~$1$ and~$2$.
    All other points in~$A$ are colored~$3$. 

    To conclude, let us show that every $J \in \mathcal{I}$ contains one color exactly once.
    By the choice of~$\mathcal{I}_M$, we know that~$J$ intersects at least one element 
    of~$\mathcal{I}_M$.
    If $I \in \mathcal{I}_M$ and $I \cap J \neq \emptyset$, then $I \cap J$ is considered in 
    $\mathcal{I}_{|I}$, and thus $J$ contains a vertex of color~$1$ or~$2$.
    In fact, $I \cap J$ contains both a vertex colored~1 and a vertex colored~2 only if it 
    contains the two inclusion-minimal intervals in $\mathcal{I}_{|I}$.
    If~$J$ moreover had another vertex colored~2, then it would also contain a vertex 
    colored~1 and, in particular, for some $I' \in \mathcal{I}_M \setminus \{I\}$, 
    we would have $I' \cap J \neq \emptyset$.
    But this would imply that~$J$ contains three inclusion-minimal intervals, contradicting
    Claim~\ref{claim:3inclusionmin}.
    Therefore, the result follows.
\end{proof}


\convex*
\begin{proof}
    The class of convex-round graphs is closed under taking induced subgraphs 
    and the bipartite convex-round graphs are exactly the biconvex graphs~\cite{BaHuYe2000}.
    The result follows from Lemmas~\ref{lem:reducetobip} 
    and~\ref{lem:biconvex}.
\end{proof}

A \defi{comparability graph} is a graph that connects pairs of elements that are 
comparable to each other in a partial order.
A partial order~$P$ has dimension two if and only if there exists a partial order~$Q$ 
on the same set of elements, such that every pair of distinct elements is comparable
in exactly one of these two partial orders.
Therefore, the permutation graphs are equivalent to the comparability graphs of
2-dimensional partially ordered sets (posets).
As another consequence of Lemma~\ref{lem:biconvex}, we show that permutation graphs 
are also $\chipcf$-bounded.

\perm*
\begin{proof}
    Bipartite permutation graphs are biconvex~\cite{SpBrSt1987}.
    This allows us to use Lemma~\ref{lem:biconvex}. Indeed we could obtain a linear $\chipcf$-bounding function using Lemma~\ref{lem:biconvex} and Lemma~\ref{lem:reducetobip}, but in order to obtain an improved 
    bound on the $\chipcf$-bounding function, we do not use Lemma~\ref{lem:reducetobip} directly.
    Suppose that~$G$ is the comparability graph of a poset~$P$.
    Take a greedy antichain partition of~$P$ by iteratively picking all the maxima. 
    It is well known (sometimes called dual Dilworth's Theorem) that the number of antichains 
    corresponds to the height $h(P)$, i.e., the length of the longest chain of~$P$.
    Let $(A_1, \ldots, A_{h(P)})$ be the antichains in our partition of~$P$, where $A_{h(P)}$ 
    contains all the maxima of~$P$.
    By Lemma~\ref{lem:biconvex} we can color each~$A_i$, with $i \in [h(P)] \setminus \{1\}$, 
    using colors from $\{3i-2,3i-1,3i\}$ in a way that every vertex of~$A_{i-1}$
    has one of these colors appearing once in its neighborhood. 
    (Note that the bipartite permutation graph induced by $A_i \cup A_{i-1}$ has no isolated 
    vertices in~$A_{i-1}$ by construction.)
    Finally, consider $G[A_{h(P)} \cup A_1]$ and note that all isolated vertices of this graph
    are isolated vertices of~$G$. 
    Hence, using Lemma~\ref{lem:biconvex} again, we can color~$A_1$ with colors from
    $\{1,2,3\}$ in a way that every non-isolated vertex of~$A_{h(P)}$ has one of 
    these colors appearing once in its neighborhood. Since we have $h(P)\le \chi(G)$, the result follows.
\end{proof}

\begin{remark}
    Note that in both Theorems~\ref{thm:convex} and~\ref{thm:perm} we also get an 
    upper bound in terms of the clique number~$\omega$.
    Namely, convex-round graphs are circular-perfect~\cite{BaHu2002} and any 
    circular-perfect graph~$G$ satisfies $\chi(G) \le \omega(G)+1$~\cite{PeWa2021}.
    Furthermore, comparability graphs are perfect, hence any permutation graph 
    satisfies $\chi(G) \le \omega(G)$.
\end{remark}

We do not know how good our upper bounds from this section are.
Let us merely establish, that on the considered classes we cannot have equality of~$\chi(G)$ and~$\chio(G)$ for all~$G$.

For each natural number $n \ge 2$ we define \defi{the cocktail party graph} $CP_{n}$ as the graph obtained from the complete graph~$K_{2n}$ on $\{0, \ldots, 2n-1\}$ with a perfect matching $M=\{(2i,2i+1)\mid i\in\{0,\ldots n-1\}$ removed.

\begin{proposition}\label{prop:lb}
    For every $n\ge 2$ the cocktail party graph $CP_n$ on~$2n$ vertices is 
    a convex-round permutation graph with $\chi(CP_n)+2 \leq \chio(CP_n)$.
\end{proposition}
\begin{proof}
    
    Ordering the vertices of $CP_n$ by their indices circularly shows that this graph is convex-round.
    The partial order $P$ given by the rule $i <_P j$ whenever $i < j-1$ and $i <_P i+1$ if~$i$ is odd
    yields a poset of width~$2$ and hence dimension~$2$,
    whose comparability graph is~$CP_n$.

    Clearly, $\chi(CP_n)=n$ and let us prove by induction on~$n$ that $\chio(CP_n) \geq n+2$.
    For $n=2$ we have that $CP_n=C_4$ for which it is known that~$4$ colors are needed 
    in an odd coloring.
    Now let $n\ge 3$ and consider~$CP_n$ with a proper odd coloring $c$.
    Remove the two independent vertices $2n-2,2n-1$ from~$CP_n$ obtaining~$CP_{n-1}$. 
    If the restriction of $c$ to~$G_{n-1}$ uses at least~$n+1$ colors, then since 
    both $2n-2,2n-1$ are connected to all vertices of~$CP_{n-1}$, the coloring $c$
     needs at least one other color, i.e., it uses at least~$n+2$ colors.

    If the restriction of $c$ to~$CP_{n-1}$ uses at most~$n$ colors, then 
    by induction hypothesis there is a vertex~$CP_{n-1}$ such that all the colors
    appearing in its neighborhood appear an even number of times.
    By symmetry, we can assume that~$0$ is such a vertex.
    Since~$N(0)$ is the union of two cliques of the same size, a set of equally colored 
    vertices in~$N(0)$ must be exactly a pair $2j,2j+1$ for $1 \le j\le n-2$.
    And on~$N(0)$ at least~$n-2$ colors are used and~$0$ uses yet another color.
    Then without loss of generality, $c(2n-2)$ appears an odd number of times in $N(0)$  and~$2n-1$ must be colored differently    from~$2n-2$.
    But none of them can use any of the already used~$n-1$ colors, hence~$n+1$
    colors are used.
    Moreover, since $W_o(c,2n-2), W_o(c,2n-1)$ are non-empty in $c$, the vertex~$1$ must use yet another color.
    Hence at least $n+2$ colors are used.
\end{proof}

\subsection{Unboundedness of $\chio$ for convex bipartite graphs}
\label{sec:biconvex}

In this subsection we prove Theorem~\ref{thm:unbounded}(a). The graphs that we construct now will be convex bipartite, hence showing that
Lemma~\ref{lem:biconvex} does not extend to this class. 
For every positive integer~$k$, let~$G_k$ be the $(A,B)$-bipartite graph where 
$A=[2^k]$, $\cI$ is the set of all non-empty intervals in~$A$, and every $I \in \cI$ 
is the neighborhood of a vertex in~$B$. 
The following corresponds to Theorem~\ref{thm:unbounded}~(a).

\begin{proposition}
\label{prop:convex}
    For every positive integer~$k$, the graph $G_k$
    \begin{enumerate}[(i)]
        \item\label{it:conv:i} is convex bipartite,
        \item\label{it:conv:ii} has $2^{k+1} + \Choose{2^k}{2}$ vertices,
        \item\label{it:conv:iv} is a grid-intersection graph,
        \item\label{it:conv:iii} is the comparability graph of a $3$-dimensional poset,
        \item\label{it:conv:v} has $\chipcf(G_k) \leq k+2$,
        \item\label{it:conv:vi} has $\chio(G_k) \geq k+2$.
    \end{enumerate}
\end{proposition}
\begin{proof}
    It is easy to see that~\eqref{it:conv:i} follows from the definition of~$G_k$.

    To see~\eqref{it:conv:ii}, note that $|A| = 2^k$ and there are $2^k + \Choose{2^k}{2}$ 
    intervals in $\mathcal{I}$.

    To see~\eqref{it:conv:iv}, simply draw every vertex of~$v \in A=[2^k]$ as a vertical ray 
    starting at the coordinate $(v,0)$ and draw every vertex of~$B$ as the horizontal interval
    $I_w \in \mathcal{I}$ corresponding to its neighborhood at a height such that 
    no two horizontal intervals intersect. 

    Now to prove \eqref{it:conv:iii}, observe that the horizontal intervals representing~$B$
    can furthermore be drawn such that if an interval is contained in another one, 
    then it is drawn lower.
    This is possible since these intervals do not intersect.
    Such graphs were studied under the name of \defi{SegRay} graphs (segmentation-ray intersection graphs),
    with a special representation for which it is known that they are comparability graphs
    of $3$-dimensional posets, see~\cite[Lemma 1]{ChFeHoWi2018} and also~\cite{ChHeOtSaUe2014}.

    For~\eqref{it:conv:v}, we prove by induction on~$k$ that $\chipcf(G_k) \leq k+2$ with the 
    additional property that only~$k+1$ colors are used on~$A$.
    This is clearly true for~$G_1$, which is a path on~5 vertices, so consider $k > 1$.
    Let $A_1 = \{1,\ldots,2^{k-1}\}$ and $A_2 = \{2^{k-1}+1,\ldots,2^k\}$, and 
    for $i=1,2$ let $B_i = \{v \in B \colon N(v) \subseteq A_i\}$ and $H_i := G_k[A_i \cup B_i]$.
    Note that each~$H_i$ is isomorphic to $G_{k-1}$ so, by inductive hypothesis, 
    let $c_i \colon V(H_i) \to [k+1]$ be a vertex coloring for~$H_i$ such that~$A_i$
    only uses colors in~$[k]$.
    Exchange the labels of colors~$k+1$ and $c_2(2^{k-1}+1)$ in~$c_2$ and observe that
    color~$k+1$ appears only once in~$A_2$.
    Now we build a coloring~$c$ for~$G_k$.
    Put $c(v) = c_i(v)$ if $v \in V(H_i)$, and put $c(v) = k+2$ if $v \notin V(H_1) \cup V(H_2)$.
    Note that all vertices in~$A$, as well as all vertices in~$B$ whose neighborhood is 
    completely in some~$A_i$, see one color exactly once, by induction.
    All other vertices in~$B$ are neighbors of vertex $2^{k-1}+1$ of~$A$, and thus they see the color~$k+1$ exactly once.

    For~\eqref{it:conv:vi}, we will show that at least~$k+1$ colors are needed by~$A$.
    Since~$B$ contains a vertex that is adjacent to all vertices of~$A$,
    this will imply the claim.
    So suppose that we have a proper odd coloring where~$A$ is colored with only~$k$ colors.
    Consider all intervals of the form $[2i-1,2i] \subseteq A$, for $i \in [2^{k-1}]$. 
    In order to have a proper odd coloring, any such interval must use at least two different 
    colors $\ci,\cb \in [k]$.
    We code this by associating to the interval $[2i-1,2i]$ the vector $x^i \in \mathbb{Z}_2^k$ 
    which is~$0$ everywhere except in the entries $\ci,\cb$, where it is~$1$.

    Note that, since we have a proper odd coloring, for each pair $j   < \ell$,  $\sum_{i=j}^{\ell} x^i \neq 0$ 
    , because otherwise the interval $[2j-1,2\ell]$ would only 
    have colors appearing at an even number of times. Moreover, this implies that for every pair $j < \ell$, we have that  $\sum_{i=j}^{2^{k-1}-1} x^i \neq \sum_{i=\ell+1}^{2^{k-1}-1}{x^i}$.
    Thus, the set $X = \{\sum_{i=j}^{2^{k-1}-1} x^i : j \in [2^{k-1}-1]\}$ has $2^{k-1}-1$ different elements.  
    On the other hand, we know $X \subseteq \{x \in \mathbb{Z}_2^k \colon \sum_{r=1}^k x_r = 0\}$, which is a subspace     of $\mathbb{Z}_2^k$ of codimension~$1$, hence has $2^{k-1}$ elements. 
    Consider now the vector $x^{2^{k-1}}$. Since the interval $[2j-1, 2^k]$ must have a color an odd number of times, then for each $j< 2^{k-1}$ we have   $\sum_{i=1}^{2^{k-1}-1} x^i \neq x^{2^{k-1}}$, which implies that vector $x^{2^{k-1}}$ does not belong to $X$. Then it must be equal to 0,  which is a contradiction to the existence of the proper odd coloring.
\end{proof}

\subsection{Separating $\chio$ and $\chipcf$ in bipartite circle graphs}
\label{sec:circle}

In this section we prove Theorem~\ref{thm:unbounded}~(b). We will construct a class of bipartite circle graphs on which the 
proper conflict-free chromatic number is unbounded, while the proper odd chromatic
number is bounded. 
To the best of our knowledge, no graphs with such properties have been shown before. 
Biconvex graphs are also a proper subclass of bipartite circle graphs~\cite{YuCh1995}
but in light of our previous results, here we will work on the difference
between these classes.
Before getting to the construction, let us give a definition and a lemma.
A family~$\cF$ of subsets of a set~$V$ is \defi{nested} if $\emptyset \notin \cF$ and 
$I \cap J \neq \emptyset \implies I \subseteq J$ or $I \supseteq J$ for all $I,J \in \cF$.

\begin{lemma}
\label{lem:nested}
    If $\cF \subseteq 2^V$ is a nested family, then $V$ can be $3$-colored such that 
    every member $I \in \cF$ contains one color an odd number of times and one color an even
    number of times.
\end{lemma}
\begin{proof}
    We order~$\cF$ by inclusion and construct a coloring bottom to top in this poset 
    by induction.
    If~$I \in \cF$ is minimal, then we can clearly $3$-color its elements such that one
    color appears an odd number of times and another one an even number of times.

    Now suppose $I \in \cF$ and the elements of~$V$ covered by the (mutually disjoint) 
    inclusion maximal intervals $I_1, \ldots, I_k$ contained in~$I$ are colored such
    that each interval contained in some $I_i$, $1 \leq i \leq k$, contains some color 
    an odd number of times.

    Permute the coloring of the elements contained in each $I_i$, $1 \leq i \leq k$,
    such that for an odd number of~$I_i$, $1 \leq i \leq k$, the number of times that color $1$ (resp color $2$) appears is odd (resp. even), and for an even number of the~$I_i$, $1 \leq i \leq k$, the number of times that color $2$ (resp color $1$) appears is odd (resp. even).
    Finally, color all elements covered by~$I$ but by none of the intervals contained 
    in~$I$ by~$3$.
    In the resulting coloring, $I$ has an odd number of elements of color~$1$ and an 
    even number of elements of color~$2$.
\end{proof}

\begin{figure}[ht!]
    \centering
    \includegraphics[width=.9\textwidth]{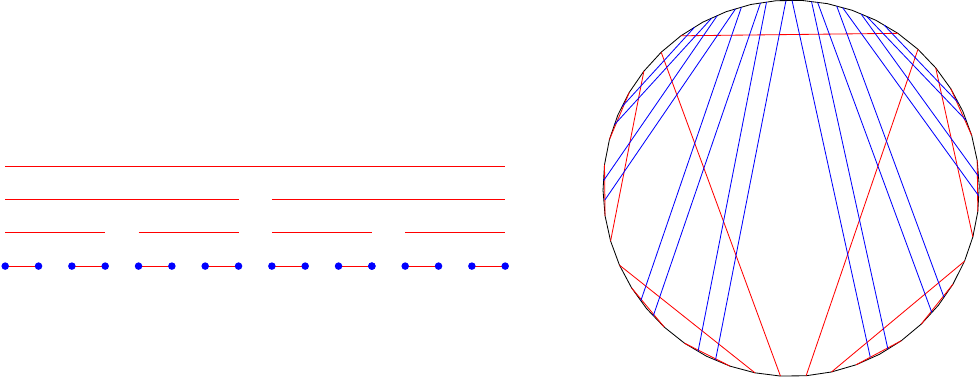}
    \caption{The interval representation for $H_4$ and its representations as chords in a circle.}
    \label{fig:bipcircle}
\end{figure}

For every positive integer~$k$, define $H_k$ as an $(A,B)$-bipartite graph where $A=[2^k]$, 
$\cI$ is the maximum nested set of intervals of length $2^i$ for all $1 \leq i \leq k$
in~$A$, i.e., greedily take intervals of the specified length starting from the beginning, 
and every $I \in \cI$ is the neighborhood of a vertex in~$B$.
See Figure~\ref{fig:bipcircle} for an illustration. 
The following corresponds to Theorem~\ref{thm:unbounded}~(b).

\begin{proposition}\label{prop:circlenotconflictfree}
    For every positive integer $k$, the graph $H_k$ satisfies the following properties:
    \begin{enumerate}[(i)]
        \item\label{it:circl:i} is a bipartite circle graph,
        \item\label{it:circl:ii} has $2^{k+1}-1$ vertices,
        \item\label{it:circl:iii} has $\chipcf(H_k) \geq k+2$,
        \item\label{it:circl:iv} has $\chio(H_k) \leq 4$.
    \end{enumerate}
\end{proposition}
\begin{proof}
    Note that~\eqref{it:circl:i} follows from the definition of~$H_k$.
    Figure~\ref{fig:bipcircle} illustrates how to represent~$H_k$ as a circle graph.

    To see~\eqref{it:circl:ii}, note that~$A$ has $2^k$ vertices and the vertices of~$B$
    corresponding to the intervals in $\cI$ ordered by inclusion form (the closure of) a full 
    binary tree with $2^{k-1}$ leaves.
    Hence, $H_k$ has $1 + 2 + 4 + \cdots + 2^k = 2^{k+1}-1$ vertices.

    To show~\eqref{it:circl:iii}, we prove by induction on~$k$ that at least~$k+1$ colors are
    needed to color~$A$ in any proper conflict-free coloring. 
    Then, since~$B$ contains a vertex incident to all vertices of~$A$, this implies that,
    in total, at least~$k+2$ colors are needed.
    
    The claim clearly holds for~$H_1$, which is a path on~$3$ vertices and needs~$2$ colors
    on its endpoints. 
    Now let $k > 1$ and suppose we have a proper conflict-free coloring of~$H_k$ such that~$A$ uses 
    only colors in~$[k]$.
    Because of the vertex in~$B$ whose neighborhood is~$A$, we can also assume, w.l.o.g,
    that color~$k$ appears exactly once in~$A$.
    Then, in one of the sets $A_1 = \{1, \ldots, 2^{k-1}\}$, $A_2 = \{2^{k-1}+1, \ldots, 2^k\}$, 
    the color~$k$ does not appear, w.l.o.g., suppose~$A_1$.
    Let $B_1 = \{v \in B \colon N(v) \subseteq A_1\}$ and $G := H_k[A_1 \cup B_1]$.
    Note that~$G$ is isomorphic to~$H_{k-1}$ and that the proper odd coloring for~$H_k$ restricted 
    to~$G$ is a proper odd coloring for~$G$ where~$A_1$ uses only colors in~$[k-1]$, which is
    a contradiction to the induction hypothesis.

    To see~\eqref{it:circl:iv}, note that, by Lemma~\ref{lem:nested}, we can color~$A$
    with~$3$ colors.
    Now for every interval~$I_w \in \cI$ of size two, color the corresponding vertex~$w \in B$
    with a color different from the two colors used in~$I_w$.
    Color all remaining vertices of~$B$ with color~$4$.
    This is a proper odd coloring for~$H_k$.
\end{proof}

While we have separated proper odd and  conflict-free chromatic numbers on 
the class of circle graphs, the following remains open.

\begin{question}\label{quest:circle}
    Are circle graphs $\chio$-bounded?
\end{question}

In order to answer this question, it is sufficient to study whether bipartite circle graphs
have bounded odd chromatic number. For this, representation of such graphs through spanning trees of planar graphs, due to de Fraysseix~\cite[Proposition 6]{deF1981}, may possibly be helpful.

\section{Acknowledgments}
The authors thank the referee for the careful reading and the suggestions that improved the presentation of the paper.

This work was initiated during the ChiPaGra workshop 2023 in Valpara\'\i so, Chile: the authors thank the organizers and the attendees for the great atmosphere.
This work is supported by the following grants and projects: 
\begin{inparaenum}[1)]
\item FAPESP-ANID Investigaci\'on Conjunta grant 2019/13364-7.
\item MATHAMSUD MATH210008. 
\item Coordena\c cão de Aperfei\c coamento de Pessoal de N\'\i vel Superior -- Brasil -- CAPES -- Finance Code 001.
\item AJ was supported by ANID/Fondecyt Regular grant 1220071 and ANID-MILENIO-NCN2024-103. 
\item KK was supported by the MICINN grant PID2022-137283NB-C22 and by the Spanish State Research Agency, through the Severo Ochoa and María de Maeztu Program for Centers and Units of Excellence in R\&D (CEX2020-001084-M).
\item CNL and YW were supported by the National Council for Scientific and Technological Development of Brazil (CNPq) Proc.~312026/2021-8 and Proc.~311892/2021-3; and joint project Proc.~404315/2023-2.
\item MS and YW were partially supported by FAPESP (Proc.~2023/03167-5).
\item DAQ and WY were supported by ANID/Fondecyt Iniciación 11201251 grant.
\item JPP was supported by ANID Doctoral Fellowship grant 21211955.
\item JZ was supported by ANID/CONICYT Fondecyt Regular grant 1241663.
\item MM  was supported by ANID Basal Grant CMM FB210005.
\item WY was supported by ANR (France) project HOSIGRA (ANR-17-CE40-0022), the National Natural Science Foundation of China (NSFC, No.~12401471) and Zhejiang Provincial Natural Science Foundation of China (ZJNSFC, No.~QN25A010023).

\end{inparaenum}

\bibliographystyle{splncs04}
\bibliography{bibliography}

\end{document}